\documentclass[10pt,a4paper,fleqn]{article}
\usepackage{amssymb,amsthm,amsfonts,hyperref}
\usepackage{amsmath}
\newtheorem{theorem}{Theorem}

\numberwithin{equation}{section}
\numberwithin{lemma}{section}
\numberwithin{theorem}{section}
\numberwithin{corollary}{section}

\allowdisplaybreaks
\begin{document}
\setcounter{page}{1}

\title{On the recursion formulas of Lauricella matrix functions}
\author{Ashish Verma\footnote{Department of Mathematics,  Prof. Rajendra Singh (Rajju Bhaiya) Institute of Physical Sciences for Study and Research,  V.B.S. Purvanchal University, Jaunpur  (U.P.)- 222003, India. \newline E-mail:	vashish.lu@gmail.com (Corresponding author)}, Ravi Dwivedi\footnote {Department of Mathematics, National Institute of Technology, Kurukshetra, India. \newline E-mail: dwivedir999@gmail.com}  
\ and Vivek Sahai\footnote {Department of Mathematics and Astronomy Lucknow University Lucknow 226007, India. \newline E-mail: sahai\_vivek@hotmail.com } }

\maketitle
\begin{abstract}
 In this paper, we find the recursion formulas for generalized Lauricella  matrix function. We also give the recursion formulas for the three variable Lauricella matrix functions.
 
 \medskip
 \noindent\textbf{Keywords}: Matrix functional calculus, Lauricella hypergeometric functions.

 \medskip
 \noindent\textbf{AMS Subject Classification}: 15A15; 33C65; 33C70. 
\end{abstract}

\section{Introduction}
The theory of matrix special functions has attracted considerable attention in the last two decades. Special matrix functions appear in the literature related to Statistics \cite{A}, Lie theory \cite{AT} and in connection with the second order matrix differential equations satisfied by matrix polynomials as well as matrix functions, for more detail see \cite{RD}--\cite{RD3}, \cite{LR}--\cite{LC1}. Recursion formulas for the Appell functions have been studied in the literature, see Saad and Srivastava \cite{OP} and Wang \cite{XW}. Authors carried out a systematic study of recursion formulas for the multivariable hypergeometric functions, \cite{VS1}-\cite{VS5}.

 Abd-Elmageed \textit{et. al.} \cite{AAAS} have obtained numerous contiguous and recursion formulas satisfied by the first Appell matrix function, namely $F_1$. Recently, Sahai and Verma obtain recursion formulas for one and two variable  hypergeometric matrix  functions in \cite{VS2}.  In the present paper, we study recursion formulas for  Lauricella  matrix function of three variables as well as $k$-variables. The paper is organized as follows.
 
  In Section~2, we give a brief review of basic definitions  that are needed in the sequel. In Section~3, we obtain the recursion formulas for generalized Lauricella  matrix function of $k$-variables. Finally, in Section~4, the recursion formulas for three variable Lauricella matrix functions are given. 
  
\section{Preliminaries}
Let the spectrum of a matrix $A$ in $\mathbb{C}^{r\times r}$, denoted by $\sigma(A)$, be the set of all eigenvalues of $A$ and let $\alpha(A) = \max \{\,\Re(z) \mid z \in \sigma(A)\,\}$ and  $\beta(A) = \min \{\,\Re(z) \mid z \in \sigma(A)\,\}$. Then for a positive stable matrix $A \in \mathbb{C}^{r \times r}$, that is $\beta(A) > 0$, the gamma matrix function is defined as \cite{LC1}
\[ \Gamma(A) = \int_{0}^{\infty} e^{-t} \, t^{A-I}\, dt
\]
and the reciprocal gamma matrix function is defined as \cite{LC}
\begin{equation}
	\Gamma^{-1}(A)= A(A+I)\dots (A+(n-1)I)\Gamma^{-1}(A+nI) , \  n\geq 1.\label{eq.07}
\end{equation}
The Pochhammer symbol  for  $A\in \mathbb{C}^{r\times r}$ is given by \cite{LC1}
\begin{equation}
	(A)_n = \begin{cases}
		I, & \text{if $n = 0$,}\\
		A(A+I) \dots (A+(n-1)I), & \text{if $n\geq 1$}.
	\end{cases}\label{c1eq.09}
\end{equation}
This gives
\begin{equation}
	(A)_n = \Gamma^{-1}(A) \ \Gamma (A+nI), \qquad n\geq 1.\label{c1eq.010}
\end{equation} 
The reciprocal gamma function $\Gamma^{-1}(z)=1/\Gamma(z)$ is an entire function of the complex variable $z$. The image of $\Gamma^{-1}(z)$ acting on $A$, denoted by $\Gamma^{-1}(A)$, is a well defined matrix. If $A+nI$ is invertible for all integers $n\geq 0$, then the reciprocal gamma function \cite{LC1} is defined as 
$\Gamma^{-1}(A)= (A)_n \ \Gamma^{-1}(A+nI)$, where $(A)_n$ is the shifted factorial matrix function for $A\in\mathbb{C}^{r\times r}$ given by \cite{LC}
\begin{align*} 
	(A)_n=
\begin{cases}
I ,
& n=0,\\
 A(A+I) \cdots (A+(n-1)I) ,
&n\geq 1 .
\end{cases}
\end{align*}
The matrix version of four Appell functions has been defined and studied in \cite{QM, RD, RD3}. The four Appell matrix functions $F_1$, $F_2$, $F_3$ and $F_4$ have been generalized to Lauricella matrix functions of $k$-variables as follows \cite{RD1, RD2}:
\begin{eqnarray}
&&\notag F_\mathcal{A}\left[{A},{B_{1}},\dots,{B_{k}};{C_{1}},\dots,{C_{k}};{x_{1}},\dots,{x_{k}}\right]\\
&&=\sum_{m_{1},\dots,{m_k}=0}^{\infty}{\left({A}\right)_{m_{1}+\dots+m_{k}}\prod_{i=1}^{k}\left({B_{i}}\right)_{m_{i}}}{\prod_{i=1}^{k}\left({C_{i}}\right)^{-1}_{m_{i}}}\prod_{i=1}^{k}\frac{x_{i}^{m_{i}}}{{m_{i}}{!}},\label{c1}
\end{eqnarray}
\begin{eqnarray}
&&\notag F_\mathcal{B}\left[{A_{1}},\dots,{A_{k}},{B_{1}},\dots,{B_{k}};{C};{x_{1}},\dots,{x_{k}}\right]\\
&&=\sum_{m_{1},\dots,{m_k}=0}^{\infty}{\prod_{i=1}^{k}\left(A_{i}\right)_{m_{i}}\prod_{i=1}^{k}\left(B_{i}\right)_{m_{i}}}{\left(C\right)^{-1}_{m_{1}+\dots+m_{k}}}\prod_{i=1}^{k}\frac{x_{i}^{m_{i}}}{{m_{i}}{!}},\label{c2}
\end{eqnarray}
\begin{eqnarray}
&&\notag F_\mathcal{C}\left[{A},{B};{C_{1}},\dots,{C_{k}};{x_{1}},\dots,{x_{k}}\right]\\
&&=\sum_{m_{1},\dots,{m_k}=0}^{\infty}{\left(A\right)_{m_{1}+\dots+m_{k}}\left({B}\right)_{m_{1}+\dots+m_{k}}}{\prod_{i=1}^{k}\left({C_{i}}\right)^{-1}_{m_{i}}}\prod_{i=1}^{k}\frac{x_{i}^{m_{i}}}{{m_{i}}{!}},\label{c3}
\end{eqnarray}
\begin{eqnarray}
&&\notag F_\mathcal{D}\left[{A},{B_{1}},\dots,{B_{k}};{C};{x_{1}},\dots,{x_{k}}\right]\\
&&=\sum_{m_{1},\dots,{m_k}=0}^\infty{\left({A}\right)_{m_{1}+\dots+m_{k}}\prod_{i=1}^{k}\left({B_{i}}\right)_{m_{i}}}{\left({C}\right)^{-1}_{m_{1}+\dots+m_{k}}}\prod_{i=1}^{k}\frac{x_{i}^{m_{i}}}{{m_{i}}{!}}.\label{c4}
\end{eqnarray}
where $A, A_1,\dots A_{k}, B, B_{1},\dots, B_{k}, C, C_{1},\dots, C_{k}$ are  matrices in $\mathbb{C}^{r\times r}$ such that $C+nI$, $C_{i}+nI$, $1 \leq i\leq k$ are invertible for all integers $n\geq 0$ and $x_1,\dots, x_{k}$ are complex variables.

The family of fourteen Lauricella  hypergeometric matrix  functions of three variables, denoted by  $F_1, \, \cdots, {F_{14}}$ ,  is given by \cite{RD1, RD2}. 
\begin{align}
&\notag F_{1}\left[{A},{B_{1}},{B_{2}},{B_{3}};{C_{1}},{C_{2}},{C_{3}};{x_{1}},{x_{2}},{x_{3}}\right]\\
&=\sum_{m_{1},m_{2},m_{3}=0}^{\infty}{\left({A}\right)_{m_{1}+m_{2}+m_{3}}\prod_{i=1}^{3}\left({B_{i}}\right)_{m_{i}}}{\prod_{i=1}^{3}\left({C_{i}}\right)^{-1}_{m_{i}}}\prod_{i=1}^{3}\frac{x_{i}^{m_{i}}}{{m_{i}}{!}};\label{c117}
\end{align}
\begin{align}
&\notag F_{2}\left[{A_{1}},{A_{2}},{A_{3}},{B_{1}},{B_{2}},{B_{3}};{C};{x_{1}},{x_{2}},{x_{3}}\right]\\
&=\sum_{m_{1},m_{2},{m_3}=0}^{\infty}{\prod_{i=1}^{3}\left(A_{i}\right)_{m_{i}}\prod_{i=1}^{3}\left(B_{i}\right)_{m_{i}}}{\left(C\right)^{-1}_{m_{1}+m_{2}+m_{3}}}\prod_{i=1}^{3}\frac{x_{i}^{m_{i}}}{{m_{i}}{!}};\label{c118}
\end{align}
\begin{align}
&\notag F_{3}\left[A_{1}, A_{2}, A_{2}, B_{1},   B_{2}, B_{1}; C_{1}, C_{2}, C_{3};{x_{1}},{x_{2}},{x_{3}}\right]\\
&=\sum_{m_{1},m_{2},m_{3}=0}^{\infty}{(A_{1})_{m_{1}}(A_{2})_{m_{2}+m_{3}}(B_{1})_{m_{1}+m_{3}}(B_{2})_{m_{2}}}{\prod_{i=1}^{3}(C_{i})^{-1}_{m_i}}\prod_{i=1}^{3}\frac{x_{i}^{m_{i}}}{{m_{i}}{!}};\label{c124}
\end{align}
\begin{align}
&\notag F_{4}\left[{A_{1}},{A_{1}},{A_{1}},{B_{1}},{B_{2}},{B_{2}};{C_{1}},{C_{2}},{C_{3}};{x_{1}},{x_{2}},{x_{3}}\right]\\
&=\sum_{m_{1},m_{2},m_{3}=0}^{\infty}{\left({A_{1}}\right)_{m_{1}+m_{2}+m_{3}}\left({B_{1}}\right)_{m_{1}}\left({B_{2}}\right)_{m_{2}+m_{3}}}{\prod_{i=1}^{3}\left({C_{i}}\right)^{-1}_{m_{i}}}\prod_{i=1}^{3}\frac{x_{i}^{m_{i}}}{{m_{i}}{!}};\label{c121}
\end{align}
\begin{align}
&\notag F_{5}\left[{A},{B};{C_{1}},{C_{2}},{C_{3}};{x_{1}},{x_{2}},{x_{3}}\right]\\
&=\sum_{m_{1},m_{2},m_{3}=0}^{\infty}{\left(A\right)_{m_{1}+m_{2}+m_{3}}\left({B}\right)_{m_{1}+m_{2}+m_{3}}}{\prod_{i=1}^{3}\left({C_{i}}\right)^{-1}_{m_{i}}}\prod_{i=1}^{3}\frac{x_{i}^{m_{i}}}{{m_{i}}{!}};\label{c119}
\end{align}
\begin{align}
&\notag F_{6}\left[A_{1}, A_{2}, A_{3}, B_{1}, B_{2}, B_{1}; C_{1}, C_{2}, C_{2};{x_{1}},{x_{2}},{x_{3}}\right]\\
&=\sum_{m_{1},m_{2},m_{3}=0}^{\infty}{\prod_{i=1}^{3}\left(A_{i}\right)_{m_{i}}\left(B_{1}\right)_{m_{1}+m_{3}}\left(B_{2}\right)_{m_{2}}}{\left(C_{1}\right)^{-1}_{m_{1}}\left(C_{2}\right)^{-1}_{m_{2}+m_{3}}}\prod_{i=1}^{3}\frac{x_{i}^{m_{i}}}{{m_{i}}{!}};\label{c126}
\end{align}
\begin{align}
&\notag F_{7}\left[A_{1}, A_{2}, A_{2}, B_{1}, B_{2}, B_{3}; C_{1}, C_{1}, C_{1};{x_{1}},{x_{2}},{x_{3}}\right]\\
&=\sum_{m_{1},m_{2},m_{3}=0}^{\infty}{\left(A_{1}\right)_{m_{1}}\left(A_{2}\right)_{m_{2}+m_{3}}\prod_{i=1}^{3}\left(B_{i}\right)_{m_{i}}}{\left(C_{1}\right)^{-1}_{m_{1}+m_{2}+m_{3}}}\prod_{i=1}^{3}\frac{x_{i}^{m_{i}}}{{m_{i}}{!}};\label{c129}
\end{align}
\begin{align}
&\notag F_{8}\left[A_{1}, A_{1}, A_{1}, B_{1}, B_{2}, B_{3}; C_{1}, C_{2}, C_{2};{x_{1}},{x_{2}},{x_{3}}\right]\\
&=\sum_{m_{1},m_{2},m_{3}=0}^{\infty}{\left(A_{1}\right)_{m_{1}+m_{2}+m_{3}}\prod_{i=1}^{3}\left(B_{i}\right)_{m_{i}}}{\left(C_{1}\right)^{-1}_{m_{1}}\left(C_{2}\right)^{-1}_{m_{2}+m_{3}}}\prod_{i=1}^{3}\frac{x_{i}^{m_{i}}}{{m_{i}}{!}};\label{c123}
\end{align}
\begin{align}
&\notag F_{9}\left[{A},{B_{1}},{B_{2}},{B_{3}};{C};{x_{1}}{x_{2}},{x_{3}}\right]\\
&=\sum_{m_{1},m_{2},m_{3}=0}^\infty{\left({A}\right)_{m_{1}+m_{2}+m_{3}}\prod_{i=1}^{3}\left({B_{i}}\right)_{m_{i}}}{\left({C}\right)^{-1}_{m_{1}+m_{2}+m_{3}}}\prod_{i=1}^{3}\frac{x_{i}^{m_{i}}}{{m_{i}}{!}};\label{c120}
\end{align}
\begin{align}
&\notag F_{10}\left[A_{1}, A_{2}, A_{1}, B_{1}, B_{2}, B_{1}; C_{1}, C_{2}, C_{2};{x_{1}},{x_{2}},{x_{3}}\right]\\
&=\sum_{m_{1},m_{2},m_{3}=0}^{\infty}{\left(A_{1}\right)_{m_{1}+m_{3}}\left(A_{2}\right)_{m_{2}}\left(B_{1}\right)_{m_{1}+m_{3}}\left(B_{2}\right)_{m_{2}}}{\left(C_{1}\right)^{-1}_{m_{1}}\left(C_{2}\right)^{-1}_{m_{2}+m_{3}}}\prod_{i=1}^{3}\frac{x_{i}^{m_{i}}}{{m_{i}}{!}};\label{c128}
\end{align}
\begin{align}
&\notag F_{11}\left[A_{1}, A_{2}, A_{2}, B_{1}, B_{2}, B_{1}; C_{1}, C_{2}, C_{2};{x_{1}},{x_{2}},{x_{3}}\right]\\
&=\sum_{m_{1},m_{2},m_{3}=0}^{\infty}{\left(A_{1}\right)_{m_{1}}\left(A_{2}\right)_{m_{2}+m_{3}}\left(B_{1}\right)_{m_{1}+m_{3}}
\left(B_{2}\right)_{m_{2}}}{\left(C_{1}\right)^{-1}_{m_{1}}\left(C_{2}\right)^{-1}_{m_{2}+m_{3}}}\prod_{i=1}^{3}\frac{x_{i}^{m_{i}}}{{m_{i}}{!}};\label{c125}
\end{align}
\begin{align}
&\notag F_{12}\left[A_{1}, A_{2}, A_{1}, B_{1}, B_{1}, B_{2}; C_{1}, C_{2}, C_{2};{x_{1}},{x_{2}},{x_{3}}\right]\\
&=\sum_{m_{1},m_{2},m_{3}=0}^{\infty}{\left(A_{1}\right)_{m_{1}+m_{3}}\left(A_{2}\right)_{m_{2}}\left(B_{1}\right)_{m_{1}+m_{2}}\left(B_{2}\right)_{m_{3}}}{\left(C_{1}\right)^{-1}_{m_{1}}\left(C_{2}\right)^{-1}_{m_{2}+m_{3}}}\prod_{i=1}^{3}\frac{x_{i}^{m_{i}}}{{m_{i}}{!}};\label{c127}
\end{align}
\begin{align}
&\notag F_{13}\left[A_{1}, A_{2}, A_{2}, B_{1}, B_{2}, B_{1}; C_{1}, C_{1}, C_{1};{x_{1}},{x_{2}},{x_{3}}\right]\\ 
&=\sum_{m_{1},m_{2},m_{3}=0}^{\infty}{\left(A_{1}\right)_{m_{1}}(A_{2})_{m_{2}+m_{3}}\left(B_{1}\right)_{m_{1}+m_{3}}\left(B_{2}\right)_{m_{2}}}{\left(C_{1}\right)^{-1}_{m_{1}+m_{2}+m_{3}}}\prod_{i=1}^{3}\frac{x_{i}^{m_{i}}}{{m_{i}}{!}};\label{c130}
\end{align} 
\begin{align}
&F_{14}\left[{A_{1}},{A_{1}},{A_{1}},{B_{1}},{B_{2}},{B_{1}};{C_{1}},{C_{2}},{C_{2}};{x_{1}},{x_{2}},{x_{3}}\right]\nonumber\\
&=\sum_{m_{1},m_{2},m_{3}=0}^{\infty}{\left({A_{1}}\right)_{m_{1}+m_{2}+m_{3}}\left({B_{1}}\right)_{m_{1}+m_{3}}\left({B_{2}}\right)_{m_{2}}}{\left({C_{1}}\right)^{-1}_{m_{1}}\left({C_{2}}\right)^{-1}_{m_{2}+m_{3}}}\prod_{i=1}^{3}\frac{x_{i}^{m_{i}}}{{m_{i}}{!}};\label{c122}
\end{align}
where  $A_{i}$, $B_{i}$ and $C_{i}$,  $1 \leq i \leq 3$ are matrices in $\mathbb{C}^{r\times r}$ such that $C+nI$, $C_{i}+nI$  are invertible for all integers $n\geq 0$ and $x_1$, $x_{2}$, $x_{3}$ are complex variables.

Following abbreviated notations are used in this paper. We, for example, write 
$F_\mathcal{A}$ for the series $F_\mathcal{A}\left[{A},{B_{1}},\dots,{B_{k}};{C_{1}},\dots,{C_{k}};{x_{1}},\dots,{x_{k}}\right]$ 
and \\
$F_\mathcal{A}[A+nI]$ for  $F_\mathcal{A}\left[{A}+nI,{B_{1}}, \dots,{B_{k}};{C_{1}},\dots,{C_{k}};{x_{1}},\dots,{x_{k}}\right]$.
The notation 
$F_\mathcal{A}[A+nI, B_1+n_1 I]$ stands for  $F_\mathcal{A}\left[A+nI,{B_{1}}+n_1 I, { B_{2}}, \dots,{ B_{k}};{C_{1}}, \dots,{C_{3}};{x_{1}},\dots,{x_{k}}\right]$ and  
$F_\mathcal{A}[A+nI, B_1+n_1 I, C_1+n_2 I]$ stands for \\ $F_\mathcal{A}\left[{A}+nI,{B_{1}}+n_1 I,{B_{2}},\dots,{b_{k}};{C_{1}}+n_2I,{C_{2}},\dots,{C_{k}};{x_{1}},\dots,{x_{k}}\right]$, etc. Throughout this  article we use the notation
$N_{j}=\sum_{i=1}^{j}n_{i}$, $N_0=0$,  where $j=1,\dots,k$.
\section{Recursion formulas for generalized Lauricella  matrix function}

In this section, we obtain the recursion formulas for  generalized Lauricella  matrix function.  Throughout this  paper $n$ denotes a non-negative integer.\\
{\bf Recursion formulas for  Lauricella  matrix function $F_{A}$:}
\begin{theorem}Let $A+nI$ be an invertible matrix for all $n\geq0$ and let $A B_i =B_i  A$, $B_iB_{j} = B_{j} B_i$ and $C_i C_j=C_j C_i$, $i, j=1,\dots,k$ then the following recursion formula holds true for generalized Lauricella  matrix function $F_\mathcal{A}$:
\begin{align}
&\notag F_\mathcal{A}\left[{A}+nI\right]\\
&=F_\mathcal{A}+ x_1 B_1 \left[\sum_{n_1=1}^{n}F_\mathcal{A}\left[{A}+n_1 I, B_1+I, C_1+I\right]\right] C^{-1}_1\notag\\&+  x_2 B_2 \left[\sum_{n_1=1}^{n}F_\mathcal{A}\left[{A}+n_1 I, B_2+I, C_2+I\right]\right]C^{-1}_{2}+\cdots\cdots\notag\\
&+x_k B_k\left[\sum_{n_1=1}^{n}F_\mathcal{A}\left[{A}+n_1 I, B_k+I, C_k+I\right]\right]C^{-1}_k,\label{c4eq1}\end{align}
Furthermore, if $A-n_1 I$ is invertible for every $n_1\leq n$, then
\begin{align}&
\notag F_\mathcal{A}\left[{A}-n I\right]\\
&=F_\mathcal{A}- x_1  B_{1}\left[\sum_{n_1=0}^{n-1}F_\mathcal{A}\left[{A}-n_1 I, B_1+I,  C_1+I\right]\right]C^{-1}_1\notag\\&- x_2 B_2\left[\sum_{n_1=0}^{n-1}F_\mathcal{A}\left[{A}-n_1I, B_2+I, C_2+I\right]\right] C^{-1}_2-\cdots\cdots\notag\\
&-x_k  B_k\left[\sum_{n_1=0}^{n-1}F_\mathcal{A}\left[{A}-n_1 I, B_k+I, C_k+I\right]\right] C^{-1}_k.\label{c43eq2}\end{align}
\end{theorem}
Proof: From the definition of  $F_\mathcal{A}$ and the relation 
\begin{align}\notag(A+I)_{m_1+\cdots+m_k}= A^{-1}(A)_{m_1+\cdots+m_k}\left(A+{m_1}I+\cdots+{m_k}I\right)\end{align}
we obtain the following contiguous relation:
\begin{align}&
\notag F_\mathcal{A}\left[{A}+I\right]\\
&= F_\mathcal{A}+x_1  B_1 F_\mathcal{A}\left[{A}+I, B_1+I, C_1+I\right]C^{-1}_1+ x_2  B_2 F_\mathcal{A}\left[{A}+I, B_2+I, C_2+I\right]C^{-1}_2\notag\\
&+\cdots+x_k  B_k F_\mathcal{A}\left[{A}+I, B_k+I, C_k+I\right]C^{-1}_k.\label{c4meq1}
 \end{align} 
To calculate contiguous relation for  $F_\mathcal{A}(A+2I)$,  we replace $A\mapsto A+I$ in (\ref{c4meq1}) and again using   (\ref{c4meq1}), we get 
 \begin{align}&
\notag F_\mathcal{A}\left[{A}+2I\right]\\
&= F_\mathcal{A}+ x_1 B_1 \Big[F_\mathcal{A}\left[{A}+I,  B_1+I, C_1+I\right]+F_\mathcal{A}\left[{A}+2I, B_1+I, C_1+I\right]\Big]C^{-1}_1\notag\\
&+x_2  B_2 \Big[F_\mathcal{A}\left[{A}+I,  B_2+I, C_2+I\right]+F_\mathcal{A}\left[{A}+2I, B_2+I, C_2+I\right]\Big]C^{-1}_2\notag\\
&+\cdots+x_k B_k \Big[F_\mathcal{A}\left[{A}+I,  B_k+I, C_k+I\right]+F_\mathcal{A}\left[{A}+2I, B_k+I, C_k+I\right]\Big]C^{-1}_k.\label{c4meq2}
 \end{align} 
 Iterating this technique $n$-times on $F_\mathcal{A}$ with matrix $
A+nI$, we obtain (\ref{c4eq1}). Next replace $A\mapsto A-I$ in (\ref{c4eq1}), we get 
 \begin{align}
&
\notag F_\mathcal{A}\left[{A}-I\right]\\
&= F_\mathcal{A}-x_1  B_1 F_\mathcal{A}\left[ B_1+I, C_1+I\right]C^{-1}_1- x_2  B_2 F_\mathcal{A}\left[ B_2+I, C_2+I\right]C^{-1}_2\notag\\
&-\cdots-x_k  B_k F_\mathcal{A}\left[ B_k+I, C_k+I\right]C^{-1}_k.
\label{c4meq3}
 \end{align} 
Apply this contiguous expression $n$-times for the function $F_\mathcal{A}$,
 we get (\ref{c43eq2}).
 \begin{theorem}Let $A+nI$ be an invertible matrix for all $n\geq 0$ and let $A B_i =B_i  A$, $B_iB_{j} = B_{j} B_i$ and $C_i C_j=C_j C_i$, $i, j=1,\dots,k$ then the following recursion formula holds true for generalized Lauricella  matrix function $F_\mathcal{A}$:
\begin{align}
\notag &F_\mathcal{A}\left[{A}+nI\right]=\sum_{N_k\leq n}^{}{n\choose n_1, n_2\dots n_k}\prod_{i=1}^{k}{(B_{i})_{n_{i}}  \,x_{i}^{n_i}}\notag\\
&\times\Big[ F_\mathcal{A}\left[{A}+N_k I, B_1+n_1 I,\ldots, B_k+n_k I, C_1+n_1 I,\ldots, C_k+n_kI\right]\Big] \prod_{i=1}^{k}{(C_{i})^{-1}_{n_{i}}},\label{c43eq3}
\end{align}
Furthermore, if $A-n_1 I$ is invertible for all $n_1\leq n$, then
\begin{align}
\notag &F_\mathcal{A}\left[{A}-nI\right]
=\sum_{N_k\leq n}^{}{n\choose n_1, n_2\dots n_k}\prod_{i=1}^{k}{(B_{i})_{n_{i}}  \,(-x_{i})^{n_i}}\notag\\\notag\\
&\times \Big[ F_\mathcal{A}\left[ B_1+n_1 I,\ldots, B_k+n_k I, C_1+n_1 I,\ldots, C_k+n_kI\right]\Big] \prod_{i=1}^{k}{(C_{i})^{-1}_{n_{i}}}.
\label{c43eq4}\end{align}
\end{theorem}
Proof: The proof of (\ref{c43eq3}) is by mathematical induction on  $n\in\mathbb{N}$. For $n=1$, (\ref{c43eq3}) is true. Assuming (\ref{c43eq3}) is true for $n=m$, that is, 
\begin{align}
\notag &F_\mathcal{A}\left[{A}+mI\right]=\sum_{N_k\leq m}^{}{m\choose n_1, n_2\dots n_k}\prod_{i=1}^{k}{(B_{i})_{n_{i}}  \,x_{i}^{n_i}}\notag\\
&\times\Big[ F_\mathcal{A}\left[{A}+N_k I, B_1+n_1 I,\ldots, B_k+n_k I, C_1+n_1 I,\ldots, C_k+n_kI\right]\Big] \prod_{i=1}^{k}{(C_{i})^{-1}_{n_{i}}},\label{c4oeq1}
\end{align}

Replacing $A\mapsto A+I$ in (\ref{c4oeq1}) and using the contiguous relation (\ref{c4meq1}) and then
simplifying,  we have
\begin{align}
& F_\mathcal{A}\left[{A}+(m+1)I\right]=\sum_{N_k\leq m}^{}{m\choose n_1, n_2\dots n_k}\prod_{i=1}^{k}{(B_{i})_{n_{i}}  \,x_{i}^{n_i}}\notag\\
 &\quad\times\Big[{F_\mathcal{A}}[A+N_k I, B_1+n_1I, \dots, B_k+n_k I, C_1+n_1I\dots, C_k+n_k I]+x_1{(B_1+n_1I)}\notag\\
 &\quad\times{F_\mathcal{A}}\left[A+N_k I +I, B_1+n_1I+I, B_2+n_2I, \dots,  B_k+n_k I, C_1+n_1I+I, C_2+n_2I,\dots, C_k+n_k I\right]\notag\\
 &\quad {(C_1+n_1I)^{-1}}+\cdots+ x_k {(B_k+n_k I)}\notag\\&\times{F_\mathcal{A}}[A+N_k I+I, B_1+k_1I, \dots, B_{k-1}+n_{k-1}I, B_k+n_kI+I, \notag\\&C_1+n_1I, \dots, C_{k-1}+n_{k-1}I, C_k+n_kI+I]{(C_k+n_k I)^{-1}}\Big]\prod_{i=1}^{k}{(C_{i})^{-1}_{n_{i}}}.
\end{align}
 \begin{align}
 & F_\mathcal{A}\left[{A}+(m+1)I\right]\notag\\
 &=\sum_{N_k\leq m}^{}{m\choose n_1, n_2\dots n_k}\prod_{i=1}^{k}{(B_{i})_{n_{i}}  \,x_{i}^{n_i}}\,\notag\\
 &\quad\times F_\mathcal{A}[A+N_k I, B_1+n_1I, \dots, B_k+n_k I, C_1+n_1I\dots, C_k+n_k I]\prod_{i=1}^{k}{(C_{i})^{-1}_{n_{i}}}\notag\\
 &\quad +\sum_{N_k\leq m+1}^{}{m\choose n_1-1, n_2\dots n_k}\prod_{i=1}^{k}{(B_{i})_{n_{i}}  \,x_{i}^{n_i}}\notag\\&F_\mathcal{A}[A+N_k I, B_1+n_1I, \dots, B_k+n_k I, C_1+n_1I\dots, C_k+n_k I]\prod_{i=1}^{k}{(C_{i})^{-1}_{n_{i}}}\notag\\
 &\quad +\cdots+\sum_{N_k\leq m+1}^{}{m\choose n_1, n_2\dots n_k-1}\prod_{i=1}^{k}{(B_{i})_{n_{i}}  \,x_{i}^{n_i}}\,\notag\\&\times F_\mathcal{A}[A+N_k I, B_1+n_1I, \dots, B_k+n_k I, C_1+n_1I\dots, C_k+n_k I]\prod_{i=1}^{k}{(C_{i})^{-1}_{n_{i}}}.
\label{c4oeq3}
 \end{align}

Using Pascal's identity in (\ref{c4oeq3}), we  have  
 \begin{align}
\notag &F_\mathcal{A}\left[{A}+(m+1)I\right]=\sum_{N_k\leq m+1}^{}{m+1\choose n_1, n_2\dots n_k}\prod_{i=1}^{k}{(B_{i})_{n_{i}}  \,x_{i}^{n_i}}\notag\\
&\times\Big[ F_\mathcal{A}\left[{A}+N_k I, B_1+n_1 I,\ldots, B_k+n_k I, C_1+n_1 I,\ldots, C_k+n_kI\right]\Big] \prod_{i=1}^{k}{(C_{i})^{-1}_{n_{i}}},\label{c4oeq4}
\end{align}
This establishes (\ref{c43eq3}) for $n=m+1$. Hence by induction result given in (\ref{c43eq3}) is true for all values of $n$. The second recursion formula (\ref{c43eq4}) can be proved in an analogous manner.
\begin{theorem} Let $B_i+nI$ be an invertible matrix for all $n\geq0$ and let  $C_i C_j=C_j C_i$, $i, j=1,\dots,k$ then the following recursion formula holds true for generalized Lauricella  matrix function $F_\mathcal{A}$:
\begin{align}
 & F_\mathcal{A}\left[{B_i}+nI\right]\notag\\&= F_\mathcal{A}+ x_i  A\left[\sum_{n_1=1}^{n}F_\mathcal{A}\left[{A}+I, B_i+n_1I,  C_i+I\right]\right]C_i^{-1},
\label{c43eq5}\end{align}
Furthermore, if $B_i-n_1 I$ is invertible for $n_1\leq n$, then
\begin{align}
& F_\mathcal{A}\left[{B_i}-nI\right]\notag\\&= F_\mathcal{A}- x_i A \left[\sum_{n_1=0}^{n-1}F_\mathcal{A}\left[{A}+I, B_i-n_1 I, C_i+I\right]\right]C_i^{-1}.
\label{c43eq6}\end{align}
\end{theorem}
\begin{theorem} Let $B_i+nI$ be an invertible matrix for all $n\geq0$ and let  $C_i C_j=C_j C_i$, $i, j=1,\dots,k$ then the following recursion formula holds true for generalized Lauricella  matrix function $F_\mathcal{A}$:
\begin{align}
& F_\mathcal{A}\left[{B_i}+nI\right]\notag\\&= \sum_{n_1=0}^{n}{n\choose n_1}{(A)_{n_1}}\,{{x_i}^{n_1}}\Big[F_\mathcal{A}\left[A+n_1I, {B_i}+n_1I, C_i+n_1I\right]\Big]{(C_i)^{-1}_{n_1}},
\label{c43eq7}\end{align}
Furthermore, if $B_i-n_1 I$ is invertible for $n_1\leq n$, then
\begin{align}
& F_\mathcal{A}\left[{B_i}-nI\right]\notag\\&= \sum_{n_1=0}^{n}{n\choose n_1}{(A)_{n_1}}\,{{(-x_i)}^{n_1}}\Big[F_\mathcal{A}\left[A+n_1I,  C_i+n_1I\right]\Big]{(C_i)^{-1}_{n_1}}.
\label{c43eq8}\end{align}
\end{theorem}
\begin{theorem} Let $C_i-nI$ be an invertible matrix for all $n\geq0$ and let $A B_i =B_i  A$, $B_i B_j= B_j B_i$ and $C_i C_j=C_j C_i$, $i, j=1,\dots,k$ then the following recursion formula holds true for generalized Lauricella  matrix function $F_\mathcal{A}$:
\begin{align}
& F_\mathcal{A}\left[{C_i}-nI\right] \nonumber\\
& = F_\mathcal{A} + x_i AB_i\Big[ \sum_{n_1=1}^{n}{ F_\mathcal{A}\left[A+I,  B_i+I, {C_i}+(2-n_1)I\right] }{(C_i-n_1I)^{-1}(C_i-(n_1-1)I)^{-1}}\Big].
\label{c43eq9}\end{align}
\end{theorem}
\subsection*{ Recursion formulas for  Lauricella  matrix function $F_{\mathcal{B}}$:}
\begin{theorem}\label{kbth1} Let $A_i+nI$ be an invertible matrix for all $n\geq0$ and let $A_i B_j =B_j  A_i$, $i, j=1,\dots,k$ then the following recursion formula holds true for generalized Lauricella  matrix function $F_{\mathcal{B}}$:
\begin{align}
& F_{\mathcal{B}}\left[A_i+nI\right]=  F_{\mathcal{B}}+x_i B_i\Big[\sum_{n_1=1}^{n}F_{\mathcal{B}}\left[A_i+n_1I, B_i+I, C+I\right]\Big]C^{-1},
\label{c43eq11}\end{align}
Furthermore, if $A_i-n_1 I$ is invertible for all $n_1\leq n$, then
\begin{align}
& F_{\mathcal{B}}\left[A_i-nI\right]=  F_{\mathcal{B}}-x_i B_i\Big[\sum_{n_1=0}^{n-1}F_{\mathcal{B}}\left[A_i-n_1I, B_i+I, C+I\right]\Big]C^{-1}.
\label{c43eq12}\end{align}
\end{theorem}
\begin{theorem}\label{kbth2} Let $A_i+nI$ be an invertible matrix for all $n\geq 0$ and let $A_i B_j =B_j  A_i$, $i, j =1,\dots,k$. Then the following recursion formula holds true for generalized Lauricella  matrix function $F_{\mathcal{B}}$:
\begin{align}
& F_{\mathcal{B}}\left[A_i+nI\right]=\sum_{n_1=0}^{n}{n\choose n_1}{(B_i)_{n_{1}}}x_{i}^{n_1}\Big[F_{\mathcal{B}}\left[A_i+n_1I, B_i+n_1I, C+n_1I\right]\Big]{(C)^{-1}_{n_{1}}}.\label{c43eq13}
\end{align}
Furthermore, if $A_i-n_1 I$ is invertible for all $n_1\leq n$, then
\begin{align}
&& F_{\mathcal{B}}\left[A_i-nI\right]=\sum_{n_1=0}^{n}{n\choose n_1}{(B_i)_{n_{1}}}(-x_{i})^{n_1}\Big[F_{\mathcal{B}}\left[ B_i+n_1I, C+n_1I\right]\Big]{(C)^{-1}_{n_{1}}}.\label{c43eq14}
\end{align}
\end{theorem}
Similarly,  the recursion formulas  $F_{\mathcal{B}}\left[{B_i}+nI\right],  F_\mathcal{B}\left[{B_i}-nI\right]$  for $F_\mathcal{B}$ can be obtain by interchange of $A_i\leftrightarrow  B_i$ in Theorem~\ref{kbth1} and Theorem~\ref{kbth2}.
\begin{theorem}Let $C-nI$ be an invertible matrix for all $n\geq 0$ and let $A_i B_j =B_j  A_i$, $i, j =1,\dots,k$. Then following recursion formula holds true for generalized Lauricella  matrix function $F_{\mathcal{B}}$:
\begin{align}
\notag& F_{\mathcal{B}}\left[C-nI\right]\\&= F_{\mathcal{B}} +x_1 A_1 B_1\Big[\sum_{n_1=1}^{n}{F_{\mathcal{B}}\left[A_1+I, B_1+I, C+(2-n_1)I\right]}{(C-n_1I)^{-1}(C-(n_1-1)I)^{-1}}\Big]\notag\\
&+x_2 A_2 B_2\Big[\sum_{n_1=1}^{n}{F_{\mathcal{B}}\left[A_2+I, B_2+I, C+(2-n_1)I\right]}{(C-n_1I)^{-1}(C-(n_1-1)I)^{-1}}\Big]+\cdots\cdots\notag\\
 &+x_k A_k B_k\Big[\sum_{n_1=1}^{n}{F_{\mathcal{B}}\left[A_k+I, B_k+I, C+(2-n_1)I\right]}{(C-n_1I)^{-1}(C-(n_1-1)I)^{-1}}\Big].\label{c43eq15}
\end{align}
\end{theorem}
Proof: From the definition of $ F_{\mathcal{B}}$ and the relation 
\begin{align*}
{(C-I)^{-1}_{m_1+m_2+\cdots+m_k}}={(C)^{-1}_{m_1+m_2+\cdots+m_k}}\left[1+{m_1}{(C-I)^{-1}}+{m_2}{(C-I)^{-1}}+\cdots+{m_k}{(C-I)^{-1}}\right],
\end{align*} 
we obtain the following contiguous relation:
\begin{align}
\notag& F_{\mathcal{B}}\left[C-nI\right]\\&= F_{\mathcal{B}} +x_1 A_1 B_1\Big[{F_{\mathcal{B}}\left[A_1+I, B_1+I, C+I\right]}{(C-I)^{-1}(C)^{-1}}\Big]\notag\\
&+x_2 A_2 B_2\Big[{F_{\mathcal{B}}\left[A_2+I, B_2+I, C+I\right]}{(C-I)^{-1}(C)^{-1}}\Big]+\cdots\cdots\notag\\
 &+x_k A_k B_k\Big[{F_{\mathcal{B}}\left[A_k+I, B_k+I, C+I\right]}{(C-I)^{-1}(C)^{-1}}\Big].\label{c4aeq15}
\end{align}
Using this contiguous relation $n$-times to the matrix function $F_{\mathcal{B}}$ with matrix $C-nI$, we get  (\ref{c43eq15}).
\subsection*{Recursion formulas for $F_{\mathcal{C}}$:}
We establish the recursion formulas about  matrix $A$. The recursion formulas about $B$ can be obtain in the similar manner.
\begin{theorem} Let $A+nI$ be an invertible matrix for all  $n\geq 0$ and let $A B =B  A$. Then the following recursion formula holds true for generalized Lauricella  matrix function $F_{\mathcal{C}}$:
\begin{align}
F_{\mathcal{C}}\left[{A}+nI\right] &=F_{\mathcal{C}}+ x_1 B\Big[ \sum_{n_1=1}^{n}F_{\mathcal{C}}\left[{A}+n_1I, B+I, C_1+I\right]\Big] C^{-1}_1\notag\\
& \quad +x_2 B\Big[ \sum_{n_1=1}^{n}F_{\mathcal{C}}\left[{A}+n_1I, B+I, C_2+I\right]\Big] C^{-1}_2\notag\\
& \quad +\cdots+x_k B\Big[ \sum_{n_1=1}^{n}F_{\mathcal{C}}\left[{A}+n_1I, B+I, C_k+I\right]\Big] C^{-1}_k.\label{c43eq16}
\end{align}
Furthermore, if $A-n_1 I$ is invertible for all $n_1 \leq n$, then
\begin{align}
F_{\mathcal{C}}\left[{A}-nI\right] &=F_{\mathcal{C}}- x_1 B\Big[ \sum_{n_1=0}^{n-1}F_{\mathcal{C}}\left[{A}-n_1I, B+I, C_1+I\right]\Big] C^{-1}_1\notag\\
& \quad -x_2 B\Big[ \sum_{n_1=0}^{n-1}F_{\mathcal{C}}\left[{A}-n_1I, B+I, C_2+I\right]\Big] C^{-1}_2\notag\\
& \quad -\cdots-x_k B\Big[ \sum_{n_1=0}^{n-1}F_{\mathcal{C}}\left[{A}-n_1I, B+I, C_k+I\right]\Big] C^{-1}_k.\label{c43eq17}
\end{align}
\end{theorem}
\begin{theorem}
 Let $A+nI$ be an invertible matrix for all $n\geq 0$ and let $A B =B  A$. Then the following recursion formula holds true for generalized Lauricella  matrix function $F_{\mathcal{C}}$:
\begin{align}
\notag F_{\mathcal{C}}\left[{A}+nI\right]&=\sum_{N_k\leq n}^{}{n\choose n_1, n_2\dots n_k}\prod_{i=1}^{k}x_{i}^{n_i} {(B)}_{N_{k}}  \notag\\
&\times\Big[ F_{\mathcal{C}}\left[{A}+N_k I, B+N_k I, C_1+n_1 I,\ldots, C_k+n_kI\right]\Big] \prod_{i=1}^{k}{(C_{i})^{-1}_{n_{i}}}.\label{c43eq18}
\end{align}
Furthermore, if $A-n_1 I$ is invertible for all $n_1\leq n$, then
\begin{align}
\notag F_{\mathcal{C}}\left[{A}-nI\right]&
=\sum_{N_k\leq n}^{}{n\choose n_1, n_2\dots n_k}\prod_{i=1}^{k}  \,(-x_{i})^{n_i}(B)_{N_k}\notag\\
&\times \Big[ F_{\mathcal{C}}\left[ B+N_k I, C_1+n_1 I,\ldots, C_k+n_kI\right]\Big] \prod_{i=1}^{k}{(C_{i})^{-1}_{n_{i}}}. \label{c43eq19}
\end{align}
\end{theorem}
Similarly,  recursion formulas  $F_{\mathcal{C}}\left[{B}+nI\right],  F_\mathcal{C}\left[{B}-nI\right]$  for $F_\mathcal{C}$ can be obtained by interchange of $A \leftrightarrow  B$  in above theorems. Finally, we present the recursion formulas of  $F_{\mathcal{C}}$ about matrix  $C_i$.
\begin{theorem}
Let $C_i-nI$ be an invertible matrix for all $n\geq0$ and let $A B =B A$ and $C_i C_j=C_j C_i$, $i, j=1,\dots,k$. Then the following recursion formula holds true for generalized Lauricella  matrix function $F_{A}$:
\begin{align}
& F_{\mathcal{C}}\left[{C_i}-nI\right]\nonumber\\
& = F_{\mathcal{C}} + x_i AB\Big[ \sum_{n_1=1}^{n}{ F_{\mathcal{C}}\left[A+I,  B+I, {C_i}+(2-n_1)I\right] }{(C_i-n_1I)^{-1}(C_i-(n_1-1)I)^{-1}}\Big].
\label{c43eq20}\end{align}
\end{theorem}

\subsection*{ Recursion formulas of $F_\mathcal{D}:$}

\begin{theorem}
Let $A+nI$ be an invertible matrix for all $n\geq 0$ and let $A B_i =B_i  A$, $i=1, 2\dots, k$.  Then the following recursion formula holds true for generalized Lauricella  matrix function $F_{\mathcal{D}}$:
\begin{align}
F_{\mathcal{D}}\left[{A}+nI\right] &=F_{\mathcal{D}}+ x_1 B_1\Big[ \sum_{n_1=1}^{n}F_{\mathcal{D}}\left[{A}+n_1I, B_1+I, C+I\right]\Big] C^{-1}\notag\\
& \quad +x_2 B_2\Big[ \sum_{n_1=1}^{n}F_{\mathcal{D}}\left[{A}+n_1I, B_2+I, C+I\right]\Big] C^{-1}\notag\\
& \quad +\cdots+x_k B_k \Big[ \sum_{n_1=1}^{n}F_{\mathcal{D}}\left[{A}+n_1I, B_k+I, C+I\right]\Big] C^{-1}.\label{c43eq21}
\end{align}
Furthermore, if $A-n_1 I$ is invertible for all $n_1\leq n$, then
\begin{align}
F_{\mathcal{D}}\left[{A}-nI\right] &=F_{\mathcal{D}}- x_1 B_1\Big[ \sum_{n_1=0}^{n-1}F_{\mathcal{D}}\left[{A}-n_1I, B_1+I, C+I\right]\Big] C^{-1}\notag\\
&\quad -x_2 B_2\Big[ \sum_{n_1=0}^{n-1}F_{\mathcal{D}}\left[{A}-n_1I, B_2+I, C+I\right]\Big] C^{-1}-\cdots\cdots\notag\\
& \quad -x_k B_k\Big[ \sum_{n_1=0}^{n-1}F_{\mathcal{D}}\left[{A}-n_1I, B_k+I, C+I\right]\Big] C^{-1}. \label{c43eq22}
\end{align}
\end{theorem}
\begin{theorem}
Let $A+nI$ be an invertible matrix for all $n\geq0$ and let $A B_i =B_i  A$, $i=1, 2,\dots,k$. Then the following recursion formula holds true for generalized Lauricella  matrix function $F_{\mathcal{D}}$:
\begin{align}
\notag F_{\mathcal{D}}\left[{A}+nI\right]&=\sum_{N_k\leq n}^{}{n\choose n_1, n_2\dots n_k}\prod_{i=1}^{k}x_{i}^{n_i} {(B_i)}_{n_{i}}  \notag\\
& \quad \times\Big[ F_{\mathcal{D}}\left[{A}+N_k I, B_1+n_1 I,\dots, B_k+n_kI,  C+N_kI\right]\Big] {(C)^{-1}_{N_{k}}}.\label{c43eq23}
\end{align}
Furthermore, if $A-n_1 I$ is invertible for all $n_1\leq n$, then
\begin{align}
\notag F_{\mathcal{D}}\left[{A}-nI\right]&
=\sum_{N_k\leq n}^{}{n\choose n_1, n_2\dots n_k}\prod_{i=1}^{k}  \,(-x_{i})^{n_i}(B_i)_{n_i}\notag\\
& \quad \times \Big[ F_{\mathcal{D}}\left[ B_1+n_1 I,\dots, B_k+n_kI,  C+N_kI\right]\Big] {(C)^{-1}_{N_{k}}}.
\label{c43eq24}\end{align}
\end{theorem}
\begin{theorem}
Let $B_i+nI$ be an invertible matrix for all $n\geq 0$. Then the following recursion formula holds true for generalized Lauricella  matrix function $F_{\mathcal{D}}$:
\begin{align}
& F_{\mathcal{D}}\left[B_i+nI\right]=  F_{\mathcal{D}}+ x_i A\Big[\sum_{n_1=1}^{n}F_{\mathcal{D}}\left[A+I, B_i+n_1I,  C+I\right]\Big]C^{-1}.\label{c43eq25}
\end{align}
Furthermore, if $B_i-n_1 I$ is invertible for all $n_1\leq n$, then
\begin{align}
& F_{\mathcal{D}}\left[B_i-nI\right]=  F_{\mathcal{D}}- x_i A\Big[\sum_{n_1=0}^{n-1}F_{\mathcal{D}}\left[A+I, B_i-n_1I,  C+I\right]\Big]C^{-1}.\label{c43eq26}
\end{align}
\end{theorem}
\begin{theorem}
Let $B_i+nI$ be an invertible matrix for all $n\geq 0$. Then the following recursion formula holds true for generalized Lauricella  matrix function $F_{\mathcal{D}}$:
\begin{align}
& F_{\mathcal{D}}\left[B_i+nI\right]=\sum_{n_1=0}^{n}{n\choose n_1}(A)_{n_{1}} x_{i}^{n_1} F_{\mathcal{D}}\left[A+n_1I, B_i+n_1I, C+n_1I\right]\Big](C)^{-1}_{n_1}.\label{c43eq27}
\end{align}
Furthermore, if $B_i-n_1 I$ is invertible for all $n_1 \leq n$, then
\begin{align}
& F_{\mathcal{D}}\left[B_i-nI\right]=\sum_{n_1=0}^{n}{n\choose n_1}(A)_{n_{1}}(- x_{i})^{n_1} F_{\mathcal{D}}\left[A+n_1I, C+n_1I\right]\Big](C)^{-1}_{n_1}.\label{c43eq28}
\end{align}
\end{theorem}
\begin{theorem}
Let $C-nI$ be an invertible matrix for all $n\geq 0$ and let $A B_i =B_i  A$, $i =1,\dots,k$. Then the following recursion formula holds true for generalized Lauricella  matrix function $F_{\mathcal{D}}$:
\begin{align}
\notag& F_{\mathcal{D}}\left[C-nI\right]\\&= F_{\mathcal{D}} +x_1 A  B_1\Big[\sum_{n_1=1}^{n}{F_{\mathcal{D}}\left[A+I, B_1+I, C+(2-n_1)I\right]}{(C-n_1I)^{-1}(C-(n_1-1)I)^{-1}}\Big]\notag\\
&+x_2 A B_2\Big[\sum_{n_1=1}^{n}{F_{\mathcal{D}}\left[A+I, B_2+I, C+(2-n_1)I\right]}{(C-n_1I)^{-1}(C-(n_1-1)I)^{-1}}\Big]+\cdots\cdots\notag\\
 &+x_k A B_k\Big[\sum_{n_1=1}^{n}{F_{\mathcal{D}}\left[A+I, B_k+I, C+(2-n_1)I\right]}{(C-n_1I)^{-1}(C-(n_1-1)I)^{-1}}\Big].\label{c43eq29}
\end{align}
\end{theorem}

\section{Recursion formulas for  Lauricella  matrix function}
There are fourteen Lauricella matrix functions of three variables denoted by $F_1,  \dots,  F_{14}$. Out of these, $F_1$, $F_2$, $F_5$ and $F_9$ are particular cases of generalized Lauricella matrix functions $F_{\mathcal{A}}$, $F_{\mathcal{B}}$, $F_{\mathcal{C}}$ and $F_{\mathcal{D}}$ respectively for $n = 3$. 
We give below the recurrence formulas of remaining ten Lauricella matrix functions, \emph{viz.}, $F_3$, $F_4$, $F_6$, $F_7$, $F_8$, $F_{10}, F_{11}, F_{12}$, $F_{13}, F_{14}$. We start by finding the recursion formulas for $F_3$.

\subsection*{\bf Recursion formulas of $F_3$:}

\begin{theorem}
Let $A_1+nI$ be an invertible matrix for all $n\geq 0$ and let $A_i B_1 =B_1  A_i$ , $i=1, 2$ and $C_1 C_j=C_j C_1$, $j= 2, 3$. Then the following recursion formula holds true for Lauricella  matrix function $F_3$:
\begin{align}
& F_{3}\left[A_1+nI\right]=  F_{3}+ x_1 B_1\Big[\sum_{n_1=1}^{n}F_{3}\left[A_1+n_1 I, B_1+ I,  C_1+I\right]\Big]C^{-1}_{1}.
\label{c43eq30}\end{align}
Furthermore, if $A_1-n_1 I$ is invertible for all $n_1\leq n$, then
\begin{align}
&F_{3}\left[A_1-nI\right]=  F_{3}- x_1 B_1\Big[\sum_{n_1=0}^{n-1}F_{3}\left[A_1-n_1 I, B_1+ I,  C_1+I\right]\Big]C^{-1}_{1}.
\label{c43eq31}\end{align}
\end{theorem}
\begin{theorem}
Let $A_1+nI$ be an invertible matrix for all $n\geq0$ and let $A_i B_1 =B_1  A_i$ , $i=1, 2$ and $C_1 C_j=C_j C_1$, $j= 2, 3$. Then the following recursion formula holds true for Lauricella  matrix function $F_3$:
\begin{align}
& F_3\left[A_1+nI\right] =\sum_{n_1=0}^{n}{n\choose n_1} (B_1)_{n_1} x_{1}^{n_1} F_{3}\left[A_1+n_1I, B_1+n_1I, C_1+n_1I\right]\Big](C_1)^{-1}_{n_1}.\label{c43eq32}
\end{align}
Furthermore, if $A_1-n_1 I$ is invertible for all $n_1\leq n$, then
\begin{align}
& F_3\left[A_1-nI\right] =\sum_{n_1=0}^{n}{n\choose n_1} (B_1)_{n_1}(- x_{1})^{n_1} F_{3}\left[ B_1+n_1I, C_1+n_1I\right]\Big](C_1)^{-1}_{n_1}.
\label{c43eq33}\end{align}
\end{theorem}
\begin{theorem}
Let $A_2+nI$ be an invertible matrix for all $n\geq0$ and let $A_i B_1 =B_1  A_i$ , $i=1, 2$ ; $B_2 C_j=C_j B_2$, $j= 1, 2, 3$; $C_2 C_3=C_3 C_2$. Then the following recursion formula holds true for  Lauricella  matrix function $F_3$:
\begin{align}
 F_{3}\left[A_2+nI\right]&=  F_{3}+ x_2  \Big[\sum_{n_1=1}^{n}F_{3}\left[A_2+n_1 I, B_2+ I,  C_2+I\right]\Big]B_2 C^{-1}_{2}\notag\\&+x_3 B_1  \Big[\sum_{n_1=1}^{n}F_{3}\left[A_2+n_1 I, B_1+ I,  C_3+I\right]\Big] C^{-1}_{3}. \label{c43eq34}
\end{align}
Furthermore, if $A_2-n_1 I$ is invertible for $n_1\leq n$, then
\begin{align}
 F_{3}\left[A_2-nI\right]&=  F_{3}- x_2  \Big[\sum_{n_1=0}^{n-1}F_{3}\left[A_2-n_1 I, B_2+ I,  C_2+I\right]\Big]B_2 C^{-1}_{2}\notag\\&-x_3 B_1  \Big[\sum_{n_1=0}^{n-1}F_{3}\left[A_2-n_1 I, B_1+ I,  C_3+I\right]\Big] C^{-1}_{3}.
\label{c43eq35}\end{align}
\end{theorem}
\begin{theorem}
Let $A_2+nI$ be an invertible matrix for all $n\geq0$ and let $A_i B_1 =B_1  A_i$ , $i=1, 2$ ; $B_2 C_j=C_j B_2$, $j= 1, 2, 3$; $C_2 C_3=C_3 C_2$. Then the following recursion formula holds true for  Lauricella  matrix function $F_3$:
\begin{align}
&F_3\left[A_2+nI\right]\nonumber\\
&=\sum_{N_2\leq n}^{}{n\choose n_1, n_2} (B_1)_{n_2} x_{2}^{n_1}x_{3}^{n_2} \Big[ F_{3}\left[A_2+N_2I, B_1+n_2I, B_2+n_1I, C_2+n_1I, C_3+n_2 I\right]\Big] \notag\\
& \quad \times (B_2)_{n_1}(C_2)^{-1}_{n_1}(C_3)^{-1}_{n_2}.
\label{c43eq36}\end{align}
Furthermore, if $A_2-n_1 I$ is invertible for all $n_1\leq n$, then
\begin{align}
F_3\left[A_2-nI\right]&=\sum_{N_2\leq n}^{}{n\choose n_1, n_2} (B_1)_{n_2} (-x_{2})^{n_1}(-x_{3})^{n_2}\notag\\&\Big[ F_{3}\left[B_1+n_2I, B_2+n_1I, C_2+n_1I, C_3+n_2 I\right]\Big](B_2)_{n_1}(C_2)^{-1}_{n_1}(C_3)^{-1}_{n_2}.
\label{c43eq37}\end{align}
\end{theorem}
\begin{theorem}
Let $B_1+nI$ be an invertible matrix for all $n\geq0$ and let $A_1 A_2 = A_2  A_1$  ; $C_1 C_j=C_j C_1$, $j= 2, 3$. Then the following recursion formula holds true for  Lauricella  matrix function $F_3$:
\begin{align}
 F_{3}\left[B_1+nI\right]&=  F_{3}+ x_1 A_1  \Big[\sum_{n_1=1}^{n}F_{3}\left[A_1+ I, B_1+n_1 I,  C_1+I\right]\Big] C^{-1}_{1}\notag\\&+x_3 A_2  \Big[\sum_{n_1=1}^{n}F_{3}\left[A_2+ I, B_1+n_1 I,  C_3+I\right]\Big] C^{-1}_{3}.
\label{c43eq38}\end{align}
Furthermore, if $B_1-n_1 I$ is invertible for all $n_1\leq n$, then
\begin{align}
  F_{3}\left[B_1-nI\right]&=  F_{3}- x_1 A_1  \Big[\sum_{n_1=0}^{n-1}F_{3}\left[A_1+ I, B_1-n_1 I,  C_1+I\right]\Big] C^{-1}_{1}\notag\\&-x_3 A_2  \Big[\sum_{n_1=0}^{n-1}F_{3}\left[A_2+ I, B_1-n_1 I,  C_3+I\right]\Big] C^{-1}_{3}.
\label{c43eq39}\end{align}
\end{theorem}
\begin{theorem}
Let $B_1+nI$ be an invertible matrix for all $n\geq0$ and let $A_1 A_2 = A_2  A_1$  ; $C_1 C_j=C_j C_1$, $j= 2, 3$. Then the following recursion formula holds true for Lauricella  matrix function $F_3$:
\begin{align}
&F_3\left[B_1+nI\right]\nonumber\\
&=\sum_{N_2\leq n}^{}{n\choose n_1, n_2} (A_1)_{n_1} (A_2)_{n_2} x_{1}^{n_1}x_{3}^{n_2}\notag\\
& \quad \times \Big[ F_{3}\left[A_1+n_1I, A_2+n_2I,  B_1+N_2I, C_1+n_1I, C_3+n_2 I\right]\Big](C_1)^{-1}_{n_1}(C_3)^{-1}_{n_2}.
\label{c43eq40}\end{align}
Furthermore, if $B_1-n_1 I$ is invertible for all $n_1\leq n$, then
\begin{align}
F_3\left[B_1-nI\right]&=\sum_{N_2\leq n}^{}{n\choose n_1, n_2} (A_1)_{n_1} (A_2)_{n_2} (-x_{1})^{n_1}(-x_{3})^{n_2}\notag\\& \quad \times \Big[ F_{3}\left[A_1+n_1I, A_2+n_2I , C_1+n_1I, C_3+n_2 I\right]\Big](C_1)^{-1}_{n_1}(C_3)^{-1}_{n_2}.
\label{c43eq41}\end{align}
\end{theorem}
\begin{theorem}
Let $B_2+nI$ be an invertible matrix for all $n\geq0$ and let $A_1 A_2 = A_2  A_1$  and $C_2 C_3=C_3 C_2$. Then the following recursion formula holds true for  Lauricella  matrix function $F_3$:
\begin{align}
& F_{3}\left[B_2+nI\right]=  F_{3}+ x_2  A_2\Big[\sum_{n_1=1}^{n}F_{3}\left[A_2+ I, B_2+n_1 I,  C_2+I\right]\Big]C^{-1}_{2}.
\label{c43eq42}\end{align}
Furthermore, if $B_2-n_1 I$ is invertible for  all $n_1\leq n$, then
\begin{align}
&F_{3}\left[B_2-nI\right]=  F_{3}- x_2 A_2\Big[\sum_{n_1=0}^{n-1}F_{3}\left[A_2+ I, B_2-n_1 I,  C_2+I\right]\Big]C^{-1}_{2}.
\label{c43eq43}\end{align}
\end{theorem}
\begin{theorem}
Let $B_2+nI$ be an invertible matrix for all $n\geq0$ and let $A_1 A_2 = A_2  A_1$  and $C_2 C_3=C_3 C_2$.   Then the following recursion formula holds true for Lauricella  matrix function $F_3$:
\begin{align}
& F_3\left[B_2+nI\right]\notag\\&=\sum_{n_1=0}^{n}{n\choose n_1} (A_2)_{n_1} x_{2}^{n_1} F_{3}\left[A_2+n_1I, B_2+n_1I, C_2+n_1I\right]\Big](C_2)^{-1}_{n_1}.
\label{c43eq44}\end{align}
Furthermore, if $B_2-n_1 I$ is invertible for all $n_1\leq n$, then
\begin{align}
& F_3\left[B_2-nI\right]\notag\\&=\sum_{n_1=0}^{n}{n\choose n_1} (A_2)_{n_1}(- x_{2})^{n_1} F_{3}\left[ A_2+n_1I, C_2+n_1I\right]\Big](C_2)^{-1}_{n_1}.
\label{c43eq45}\end{align}
\end{theorem}
\begin{theorem} 
Let $C_1-nI$   be an invertible matrix for all $n\geq 0$ and let $A_i B_1 =B_1  A_i$, $C_1 C_i= C_i C_1$, $i=1,2$. Then the following recursion formula holds true for Lauricella  matrix function $F_{3}$:
\begin{align}
& F_{3}\left[{C_1}-nI\right]\nonumber\\
& = F_{3} + x_1  A_1 B_1\nonumber\\
 & \quad \times \Big[ \sum_{n_1=1}^{n}{ F_{3}\left[A_1+I,  B_1+I, {C_1}+(2-n_1)I\right] }{(C_1-n_1I)^{-1}(C_1-(n_1-1)I)^{-1}}\Big].
\label{c43eq46}\end{align}
\end{theorem}
\begin{theorem} 
Let $C_2-nI$   be an invertible matrix for all $n\geq0$ and let $A_1 A_2 = A_2  A_1$;  $B_2 C_i= C_i B_2$, $i=1, 2, 3$. Then the following recursion formula holds true for  Lauricella  matrix function $F_{3}$:
\begin{align}
& F_{3}\left[{C_2}-nI\right]\nonumber\\
& = F_{3} + x_2  A_2\nonumber\\
& \quad \times \Big[ \sum_{n_1=1}^{n}{ F_{3}\left[A_2+I,  B_2+I, {C_2}+(2-n_1)I\right] }B_2{(C_2-n_1I)^{-1}(C_2-(n_1-1)I)^{-1}}\Big].
\label{c43eq47}\end{align}
\end{theorem}
\begin{theorem}
 Let $C_3-nI$   be an invertible matrix for all $n\geq0$ and let $A_1 A_2 = A_2  A_1$; $A_i B_1= B_1 A_i$, $i=1, 2$. Then the following recursion formula holds true for  Lauricella  matrix function $F_{3}$:
\begin{align}
& F_{3}\left[{C_3}-nI\right]\nonumber\\
& = F_{3} + x_3  A_2  B_1\nonumber\\
& \quad \times \Big[ \sum_{n_1=1}^{n}{ F_{3}\left[A_2+I,  B_1+I, {C_3}+(2-n_1)I\right] }{(C_3-n_1I)^{-1}(C_3-(n_1-1)I)^{-1}}\Big].
\label{c43eq48}\end{align}
\end{theorem}

\subsection*{\bf Recursion formulas of $F_4$:}

\begin{theorem}
Let $A_1+nI$ be an invertible matrix for all $n\geq0$ and let $A_1 B_i =B_i  A_1$, $i=1, 2$;  $B_1 B_2= B_2 B_1$ and $C_i C_j=C_j C_i$, $i, j=1, 2,3$. Then the following recursion formula holds true for  Lauricella  matrix function $F_4$:
\begin{align}
\notag F_4\left[{A_1}+nI\right]& =F_4+ x_1 B_1 \left[\sum_{n_1=1}^{n}F_4\left[{A_1}+n_1 I, B_1+I, C_1+I\right]\right] C^{-1}_1\notag\\& \quad  +  x_2 B_2 \left[\sum_{n_1=1}^{n}F_4\left[{A_1}+n_1 I, B_2+I, C_2+I\right]\right]C^{-1}_{2}\notag\\
& \quad  +x_3 B_2\left[\sum_{n_1=1}^{n}F_4\left[{A_1}+n_1 I, B_2+I, C_3+I\right]\right]C^{-1}_3.\label{c43eq1}
\end{align}
Furthermore, if $A_1-n_1 I$ is invertible for all $n_1\leq n$, then
\begin{align}
  F_4\left[{A_1}-n I\right]  &= F_4- x_1  B_{1}\left[\sum_{n_1=0}^{n-1}F_4\left[{A_1}-n_1 I, B_1+I,  C_1+I\right]\right]C^{-1}_1\notag\\& \quad - x_2 B_2\left[\sum_{n_1=0}^{n-1}F_4\left[{A_1}-n_1I, B_2+I, C_2+I\right]\right] C^{-1}_2\notag\\
&\quad -x_3  B_2\left[\sum_{n_1=0}^{n-1}F_4\left[{A_1}-n_1 I, B_2+I, C_3+I\right]\right] C^{-1}_3.\label{c43eq49}\end{align}
\end{theorem}
\begin{theorem}
Let $A_1+nI$ be an invertible matrix for all $n\geq0$ and let $A_1 B_i =B_i  A_1$, $i=1, 2$;  $B_1 B_2= B_2 B_1$ and $C_i C_j=C_j C_i$, $i, j=1, 2,3$. Then the following recursion formula holds true for  Lauricella  matrix function $F_4$:
\begin{align}
F_4\left[A_1+nI\right]&=\sum_{N_3\leq n}^{}{n\choose n_1, n_2, n_3} (B_1)_{n_1} (B_2)_{n_2+n_3}\, x_{1}^{n_1}x_{2}^{n_2}x_{3}^{n_3}\notag\\& \quad \times\Big[ F_{4}\left[A_1+N_2I, B_1+n_1I , B_2+(n_2+n_3)I, C_1+n_1I, C_2+n_2 I,  C_3+n_3I\right]\Big]\notag\\& \quad \times (C_1)^{-1}_{n_1}(C_2)^{-1}_{n_2}(C_3)^{-1}_{n_3}.
\label{c43eq50}\end{align}
Furthermore, if $A_1-n_1 I$ is invertible for all $n_1\leq n$, then
\begin{align}
F_4\left[A_1-nI\right]&=\sum_{N_3\leq n}^{}{n\choose n_1, n_2, n_3} (B_1)_{n_1} (B_2)_{n_2+n_3} (-x_{1})^{n_1}(-x_{2})^{n_2}(-x_3)^{n_3}\notag\\& \quad \times \Big[ F_{4}\left[ B_1+n_1I , B_2+(n_2+n_3)I,  C_1+n_1I, C_2+n_2 I, C_3+n_3I\right]\Big]\notag\\& \quad \times (C_1)^{-1}_{n_1}(C_2)^{-1}_{n_2}(C_3)^{-1}_{n_3}.
\label{c43eq51}\end{align}
\end{theorem}
\begin{theorem}
Let $B_1+nI$ be an invertible matrix for all $n\geq 0$ and let  $C_1 C_j=C_j C_1$, $j= 2, 3$. Then the following recursion formula holds true for  Lauricella  matrix function $F_4$:
\begin{align}
& F_{4}\left[B_1+nI\right] =  F_{4}+ x_1 A_1\Big[\sum_{n_1=1}^{n}F_{4}\left[A_1+ I, B_1+ n_1I,  C_1+I\right]\Big]C^{-1}_{1}.
\label{c43eq52}\end{align}
Furthermore, if $B_1-n_1 I$ is invertible for all $n_1\leq n$, then
\begin{align}
&F_{4}\left[B_1-nI\right]=  F_{4}- x_1 A_1\Big[\sum_{n_1=0}^{n-1}F_{4}\left[A_1+  I, B_1-n_1 I,  C_1+I\right]\Big]C^{-1}_{1}.
\label{c43eq53}\end{align}
\end{theorem}
\begin{theorem}
Let $B_1+nI$ be an invertible matrix for all $n\geq 0$ and let  $C_1 C_j=C_j C_1$, $j= 2, 3$ Then the following recursion formula holds true for  Lauricella  matrix function $F_4$:
\begin{align}
& F_4\left[B_1+nI\right] =\sum_{n_1=0}^{n}{n\choose n_1} (A_1)_{n_1} x_{1}^{n_1} F_{4}\left[A_1+n_1I, B_1+n_1I, C_1+n_1I\right]\Big](C_1)^{-1}_{n_1}.
\label{c43eq54}\end{align}
Furthermore, if $B_1-n_1 I$ is invertible for all $n_1\leq n$, then
\begin{align}
& F_4\left[B_1-nI\right] = \sum_{n_1=0}^{n}{n\choose n_1} (A_1)_{n_1}(- x_{1})^{n_1} F_{4}\left[ A_1+n_1I, C_1+n_1I\right]\Big](C_1)^{-1}_{n_1}.
\label{c43eq55}\end{align}
\end{theorem}
\begin{theorem}
Let $B_2+nI$ be an invertible matrix for all $n\geq 0$ and let  $C_2 C_3=C_3 C_2$,  then the following recursion formula holds true for  Lauricella  matrix function $F_4$:
\begin{align}
 F_{4}\left[B_2+nI\right]&=  F_{4}+ x_2 A_1  \Big[\sum_{n_1=1}^{n}F_{4}\left[A_1+ I, B_2+n_1 I,  C_2+I\right]\Big] C^{-1}_{2}\notag\\& \quad +x_3 A_1  \Big[\sum_{n_1=1}^{n}F_{4}\left[A_1+ I, B_2+n_1 I,  C_3+I\right]\Big] C^{-1}_{3}.
\label{c43eq56}\end{align}
Furthermore, if $B_2-n_1 I$ is invertible for all $n_1\leq n$, then
\begin{align}
  F_{4}\left[B_2-nI\right]&=  F_{4}- x_2 A_1  \Big[\sum_{n_1=0}^{n-1}F_{4}\left[A_1+ I, B_2-n_1 I,  C_2+I\right]\Big] C^{-1}_{2}\notag\\& \quad -x_3 A_1  \Big[\sum_{n_1=0}^{n-1}F_{4}\left[A_1+ I, B_2-n_1 I,  C_3+I\right]\Big] C^{-1}_{3}.
\label{c43eq57}\end{align}
\end{theorem}
\begin{theorem}
Let $B_2+nI$ be an invertible matrix for all $n\geq 0$ and let  $C_2 C_3=C_3 C_2$,  then the following recursion formula holds true for  Lauricella  matrix function $F_4$:
\begin{align}
F_4\left[B_2+nI\right]& =\sum_{N_2\leq n}^{}{n\choose n_1, n_2} (A_1)_{N_2} x_{2}^{n_1}x_{3}^{n_2}\notag\\& \quad \times \Big[ F_{4}\left[A_1+N_2I,  B_2+N_2I, C_2+n_1I, C_3+n_2 I\right]\Big](C_2)^{-1}_{n_1}(C_3)^{-1}_{n_2}.
\label{c43eq58}\end{align}
Furthermore, if $B_2-n_1 I$ is invertible for all $n_1\leq n$, then
\begin{align}
F_4\left[B_2-nI\right]&=\sum_{N_2\leq n}^{}{n\choose n_1, n_2}  (A_1)_{N_2} (-x_{2})^{n_1}(-x_{3})^{n_2}\notag\\& \quad \times \Big[ F_{4}\left[A_1+N_2I, C_2+n_1I, C_3+n_2 I\right]\Big](C_2)^{-1}_{n_1}(C_3)^{-1}_{n_2}.
\label{c43eq59}\end{align}
\end{theorem}
\begin{theorem} 
Let $C_1-nI$   be an invertible matrix for all $n\geq0$ and let $A_1 B_1 =B_1  A_1$, $C_1C_j= C_j C_1$, $j=2, 3$. Then the following recursion formula holds true for Lauricella  matrix function $F_{4}$:
\begin{align}
& F_{4}\left[{C_1}-nI\right] \notag\\&= F_{4} + x_1  A_1 B_1 \notag\\& \quad \times\Big[ \sum_{n_1=1}^{n}{ F_{4}\left[A_1+I,  B_1+I, {C_1}+(2-n_1)I\right] }{(C_1-n_1I)^{-1}(C_1-(n_1-1)I)^{-1}}\Big].
\label{c43eq60}\end{align}
\end{theorem}
\begin{theorem}
 Let $C_2-nI$   be an invertible matrix for all $n\geq 0$ and let $B_2 C_i=C_i B_2$, $i=1,2,3$, $C_2 C_3= C_3 C_2$. Then the following recursion formula holds true for Lauricella  matrix function $F_{4}$:
\begin{align}
& F_{4}\left[{C_2}-nI\right] \notag\\&= F_{4} + x_2  A_1 \notag\\& \quad \times \Big[ \sum_{n_1=1}^{n}{ F_{4}\left[A_1+I,  B_2+I, {C_2}+(2-n_1)I\right] }B_2{(C_2-n_1I)^{-1}(C_2-(n_1-1)I)^{-1}}\Big].
\label{c43eq61}\end{align}
\end{theorem}
\begin{theorem}
 Let $C_3-nI$   be an invertible matrix for all $n\geq0$ and let $B_2 C_i=C_i B_2$, $i=1,2,3$. Then the following recursion formula holds true for Lauricella  matrix function $F_{4}$:
\begin{align}
& F_{4}\left[{C_3}-nI\right]\notag\\& = F_{4} + x_3  A_1 \notag\\& \quad \times \Big[ \sum_{n_1=1}^{n}{ F_{4}\left[A_1+I,  B_2+I, {C_3}+(2-n_1)I\right] }B_2 {(C_3-n_1I)^{-1}(C_3-(n_1-1)I)^{-1}}\Big].
\label{c43eq62}\end{align}
\end{theorem}

\subsection*{ Recursion formulas of $F_6$:}

\begin{theorem}
Let $A_1+nI$ be an invertible matrix for all $n\geq 0$ and let $A_i B_1 =B_1  A_i$, $i=1, 2,3 $; $C_1 C_2=C_2 C_1$. Then the following recursion formula holds true for Lauricella  matrix function $F_6$:
\begin{align}
& F_{6}\left[A_1+nI\right] =  F_{6}+ x_1 B_1\Big[\sum_{n_1=1}^{n}F_{6}\left[A_1+n_1 I, B_1+ I,  C_1+I\right]\Big]C^{-1}_{1}.
\label{c43eq63}\end{align}
Furthermore, if $A_1-n_1 I$ is invertible for all $n_1\leq n$, then
\begin{align}
&F_{6}\left[A_1-nI\right]=  F_{6}- x_1 B_1\Big[\sum_{n_1=0}^{n-1}F_{6}\left[A_1-n_1 I, B_1+ I,  C_1+I\right]\Big]C^{-1}_{1}.
\label{c43eq64}\end{align}
\end{theorem}
\begin{theorem}
Let $A_1+nI$ be an invertible matrix for all $n\geq0$ and let $A_i B_1 =B_1  A_i$, $i=1, 2,3 $; $C_1 C_2=C_2 C_1$. Then the following recursion formula holds true for Lauricella  matrix function $F_6$:
\begin{align}
& F_6\left[A_1+nI\right] =\sum_{n_1=0}^{n}{n\choose n_1} (B_1)_{n_1} x_{1}^{n_1} F_{6}\left[A_1+n_1I, B_1+n_1I, C_1+n_1I\right]\Big](C_1)^{-1}_{n_1}.
\label{c43eq65}\end{align}
Furthermore, if $A_1-n_1 I$ is invertible for all $n_1\leq n$, then
\begin{align}
& F_6\left[A_1-nI\right] =\sum_{n_1=0}^{n}{n\choose n_1} (B_1)_{n_1}(- x_{1})^{n_1} F_{6}\left[ B_1+n_1I, C_1+n_1I\right]\Big](C_1)^{-1}_{n_1}.
\label{c43eq66}\end{align}
\end{theorem}
\begin{theorem}
Let $A_2+nI$ be an invertible matrix for all $n\geq0$ and let $C_i B_2 =B_2  C_i$, $i=1, 2 $. Then the following recursion formula holds true for Lauricella  matrix function $F_6$:
\begin{align}
& F_{6}\left[A_2+nI\right] =  F_{6}+ x_2 \Big[\sum_{n_1=1}^{n}F_{6}\left[A_2+n_1 I, B_2+ I,  C_2+I\right]\Big]B_2C^{-1}_{2}.
\label{c43eq67}\end{align}
Furthermore, if $A_2-n_1 I$ is invertible for all $n_1\leq n$, then
\begin{align}
&F_{6}\left[A_2-nI\right]=  F_{6}- x_2 \Big[\sum_{n_1=0}^{n-1}F_{6}\left[A_2-n_1 I,  B_2+ I,  C_2+I\right]\Big]B_2C^{-1}_{2}.
\label{c43eq68}\end{align}
\end{theorem}
\begin{theorem}
Let $A_2+nI$ be an invertible matrix for all $n\geq 0$ and let $C_i B_2 = B_2  C_i$, $i=1, 2$. Then the following recursion formula holds true for Lauricella  matrix function $F_6$:
\begin{align}
& F_6\left[A_2+nI\right] =\sum_{n_1=0}^{n}{n\choose n_1}  x_{2}^{n_1} F_{6}\left[A_2+n_1I, B_2+n_1I, C_2+n_1I\right]\Big](B_2)_{n_1}(C_1)^{-1}_{n_1}.
\label{c43eq69}\end{align}
Furthermore, if $A_2-n_1 I$ is invertible for all $n_1\leq n$, then
\begin{align}
& F_6\left[A_2-nI\right] =\sum_{n_1=0}^{n}{n\choose n_1} (- x_{2})^{n_1} F_{6}\left[ B_2+n_1I, C_2+n_1I\right]\Big](B_2)_{n_1}(C_1)^{-1}_{n_1}.
\label{c43eq70}\end{align}
\end{theorem}
\begin{theorem}
Let $A_3+nI$ be an invertible matrix for all $n\geq 0$ and let $A_i B_1 =B_1  A_i$, $i=1, 2, 3 $. Then the following recursion formula holds true for Lauricella  matrix function $F_6$:
\begin{align}
& F_{6}\left[A_3+nI\right] =  F_{6} + x_3 B_1 \Big[\sum_{n_1=1}^{n}F_{6}\left[A_3+n_1 I, B_1+ I,  C_2+I\right]\Big]C^{-1}_{2}.
\label{c43eq71}
\end{align}
Furthermore, if $A_3-n_1 I$ is invertible for all $n_1\leq n$, then
\begin{align}
&F_{6}\left[A_3-nI\right]=  F_{6}- x_3 B_1 \Big[\sum_{n_1=0}^{n-1}F_{6}\left[A_3-n_1 I,  B_1+ I,  C_2+I\right]\Big]C^{-1}_{2}.
\label{c43eq72}\end{align}
\end{theorem}
\begin{theorem}
Let $A_3+nI$ be an invertible matrix for all $n\geq 0$ and let $A_i B_1 = B_1  A_i$, $i=1, 2, 3 $. Then the following recursion formula holds true for Lauricella  matrix function $F_6$:
\begin{align}
& F_6\left[A_3+nI\right] =\sum_{n_1=0}^{n}{n\choose n_1}  x_{3}^{n_1}(B_1)_{n_1} F_{6}\left[A_3+n_1I, B_1+n_1I, C_2+n_1I\right]\Big](C_2)^{-1}_{n_1}.
\label{c43eq73}
\end{align}
Furthermore, if $A_3-n_1 I$ is invertible for all $n_1\leq n$, then
\begin{align}
& F_6\left[A_3-nI\right] =\sum_{n_1=0}^{n}{n\choose n_1} (- x_{3})^{n_1}(B_1)_{n_1} F_{6}\left[ B_1+n_1I, C_2+n_1I\right]\Big](C_2)^{-1}_{n_1}.
\label{c43eq74}
\end{align}
\end{theorem}
\begin{theorem}
Let $B_1+nI$ be an invertible matrix for all $n\geq 0$ and let $A_i A_3 = A_3 A_i$, $i=1, 2$;  $C_1 C_2=C_2 C_1$.  Then the following recursion formula holds true for Lauricella  matrix function $F_6$:
\begin{align}
 F_{6}\left[B_1+nI\right]&=  F_{6}+ x_1 A_1  \Big[\sum_{n_1=1}^{n}F_{6}\left[A_1+ I, B_1+n_1 I,  C_1+I\right]\Big] C^{-1}_{1}\notag\\& \quad +x_3 A_3  \Big[\sum_{n_1=1}^{n}F_{6}\left[A_3+ I,  B_1+n_1 I,  C_2+I\right]\Big] C^{-1}_{2}
\label{c43eq75}\end{align}
Furthermore, if $B_1-n_1 I$ is invertible for all $n_1\leq n$, then
\begin{align}
  F_{6}\left[B_1-nI\right]&=  F_{6}- x_1 A_1  \Big[\sum_{n_1=0}^{n-1}F_{6}\left[A_1+ I, B_1-n_1 I,  C_1+I\right]\Big] C^{-1}_{1}\notag\\& \quad -x_3 A_3  \Big[\sum_{n_1=0}^{n-1}F_{6}\left[A_3+ I, B_1-n_1 I,  C_2+I\right]\Big] C^{-1}_{2}.
\label{c43eq76}
\end{align}
\end{theorem}
\begin{theorem}
Let $B_1+nI$ be an invertible matrix for all $n\geq 0$ and let $A_i A_3= A_3 A_i$, $i=1, 2$;  $C_1 C_2 = C_2 C_1$,  then the following recursion formula holds true for Lauricella  matrix function $F_6$:
\begin{align}
F_6\left[B_1+nI\right]&=\sum_{N_2\leq n}^{}{n\choose n_1, n_2} (A_1)_{n_1} (A_3)_{n_2} x_{1}^{n_1}x_{3}^{n_2}\notag\\& \quad  \times \Big[ F_{6}\left[A_1+n_1I, A_3+n_2I,   B_1+N_2I, C_1+n_1I, C_2+n_2 I\right]\Big](C_1)^{-1}_{n_1}(C_2)^{-1}_{n_2}.
\label{c43eq77}\end{align}
Furthermore, if $B_1-n_1 I$ is invertible for all $n_1\leq n$, then
\begin{align}
F_6\left[B_1-nI\right]&=\sum_{N_2\leq n}^{}{n\choose n_1, n_2}  (A_1)_{n_1} (A_3)_{n_2} (-x_{1})^{n_1}(-x_{3})^{n_2}\notag\\& \quad \times \Big[ F_{6}\left[A_1+n_1I, A_3+n_2I,   C_1+n_1I, C_2+n_2 I\right]\Big](C_1)^{-1}_{n_1}(C_2)^{-1}_{n_2}.
\label{c43eq78}\end{align}
\end{theorem}
\begin{theorem}
Let $B_2+nI$ be an invertible matrix for all $n\geq 0$ and let $A_1 A_2 =A_2  A_1$, then the following recursion formula holds true for  Lauricella  matrix function $F_6$:
\begin{align}
& F_{6}\left[B_2+nI\right] =  F_{6}+ x_2  A_2 \Big[\sum_{n_1=1}^{n}F_{6}\left[A_2+ I, B_2+n_1 I,  C_2+I\right]\Big]C^{-1}_{2}.
\label{c43eq79}\end{align}
Furthermore, if $B_2-n_1 I$ is invertible for all $n_1\leq n$, then
\begin{align}
&F_{6}\left[B_2-nI\right]=  F_{6}- x_2 A_2 \Big[\sum_{n_1=0}^{n-1}F_{6}\left[A_2+ I,  B_2- n_1I,  C_2+I\right]\Big]C^{-1}_{2}.
\label{c43eq80}\end{align}
\end{theorem}
\begin{theorem}
Let $B_2+nI$ be an invertible matrix for all $n\geq 0$ and let  $A_1 A_2 =A_2  A_1$, then the following recursion formula holds true for Lauricella  matrix function $F_6$:
\begin{align}
& F_6\left[B_2+nI\right] =\sum_{n_1=0}^{n}{n\choose n_1}  x_{2}^{n_1}(A_2)_{n_1} F_{6}\left[A_2+n_1I, B_2+n_1I, C_2+n_1I\right]\Big](C_2)^{-1}_{n_1}.
\label{c43eq81}\end{align}
Furthermore, if $B_2-n_1 I$ is invertible for all $n_1\leq n$, then
\begin{align}
& F_6\left[B_2-nI\right] =\sum_{n_1=0}^{n}{n\choose n_1} (- x_{2})^{n_1}(A_2)_{n_1} F_{6}\left[ A_2+n_1I, C_2+n_1I\right]\Big](C_2)^{-1}_{n_1}.
\label{c43eq82}
\end{align}
\end{theorem}
\begin{theorem}
 Let $C_1-nI$   be an invertible matrix for all $n\geq 0$ and let $A_i B_1 =B_1  A_i$, $i=1, 2, 3$, $C_1 C_2= C_2 C_1$. Then the following recursion formula holds true for Lauricella  matrix function $F_{6}$:
\begin{align}
& F_{6}\left[{C_1}-nI\right] \notag\\& = F_{6} + x_1  A_1 B_1 \notag\\& \quad \times \Big[ \sum_{n_1=1}^{n}{ F_{6}\left[A_1+I,  B_1+I, {C_1}+(2-n_1)I\right] }{(C_1-n_1I)^{-1}(C_1-(n_1-1)I)^{-1}}\Big].
\label{c43eq83}\end{align}
\end{theorem}
\begin{theorem}
Let $C_2-nI$ be an invertible matrix for all $n\geq0$ and let $A_1 A_2= A_2 A_1$; $A_i A_3= A_3 A_i$, $i=1, 2$, $B_2 C_i= C_i B_2$,  $A_j B_1= B_1 A_j$,  $i=1,2, j=1,2,3$. Then the following recursion formula holds true for Lauricella  matrix function $F_{6}$:
\begin{align}
\notag& F_{6}\left[C_2-nI\right]\\&= F_{6} +x_2 A_2 \Big[\sum_{n_1=1}^{n}{F_6\left[A_2+I, B_2+I, C_2+(2-n_1)I\right]}B_2{(C_2-n_1I)^{-1}(C_2-(n_1-1)I)^{-1}}\Big]\notag\\
&+x_3 A_3 B_1\Big[\sum_{n_1=1}^{n}{F_6 \left[A_3+I, B_1+I, C_2+(2-n_1)I\right]}{(C_2-n_1I)^{-1}(C_2-(n_1-1)I)^{-1}}\Big].\label{c43eq84}
\end{align}
\end{theorem}

\subsection*{ Recursion formulas of $F_7$:}

\begin{theorem}
Let $A_1+nI$ be an invertible matrix for all $n\geq 0$ and let $A_i B_1 =B_1  A_i$, $i=1, 2 $. Then the following recursion formula holds true for Lauricella  matrix function $F_7$:
\begin{align}
& F_{7}\left[A_1+nI\right] =  F_{7}+ x_1 B_1\Big[\sum_{n_1=1}^{n}F_{7}\left[A_1+n_1 I, B_1+ I,  C_1+I\right]\Big]C^{-1}_{1}.
\label{c43eq85}\end{align}
Furthermore, if $A_1-n_1 I$ is invertible for all $n_1\leq n$, then
\begin{align}
&F_{7}\left[A_1-nI\right]=  F_{7}- x_1 B_1\Big[\sum_{n_1=0}^{n-1}F_{7}\left[A_1-n_1 I, B_1+ I,  C_1+I\right]\Big]C^{-1}_{1}.
\label{c43eq86}\end{align}
\end{theorem}
\begin{theorem}
Let $A_1+nI$ be an invertible matrix for all $n\geq 0$ and let $A_i B_1 =B_1  A_i$, $i=1, 2 $.  Then the following recursion formula holds true for Lauricella  matrix function $F_7$:
\begin{align}
& F_7\left[A_1+nI\right] =\sum_{n_1=0}^{n}{n\choose n_1} (B_1)_{n_1} x_{1}^{n_1} F_{7}\left[A_1+n_1I, B_1+n_1I, C_1+n_1I\right]\Big](C_1)^{-1}_{n_1}.
\label{c43eq87}\end{align}
Furthermore, if $A_1-n_1 I$ is invertible for all $n_1\leq n$, then
\begin{align}
& F_7\left[A_1-nI\right] =\sum_{n_1=0}^{n}{n\choose n_1} (B_1)_{n_1}(- x_{1})^{n_1} F_{7}\left[ B_1+n_1I, C_1+n_1I\right]\Big](C_1)^{-1}_{n_1}.
\label{c43eq88}
\end{align}
\end{theorem}
\begin{theorem}
Let $A_2+nI$ be an invertible matrix for all $n\geq0$ and let $A_i B_2= B_2 A_i$, $i=1, 2$; $B_1 B_2= B_2 B_1$;  $C_1 B_3= B_3 C_1$.  Then the following recursion formula holds true for  Lauricella  matrix function $F_7$:
\begin{align}
 F_{7}\left[A_2+nI\right]&=  F_{7}+ x_2  B_2  \Big[\sum_{n_1=1}^{n}F_{7}\left[A_2+ n_1I, B_2+ I,  C_1+I\right]\Big] C^{-1}_{1}\notag\\& \quad +x_3   \Big[\sum_{n_1=1}^{n}F_{7}\left[  A_2+n_1 I, B_3+ I, C_1+I\right]\Big]B_3 C^{-1}_{1}.
\label{c43eq89}\end{align}
Furthermore, if $A_2-n_1 I$ is invertible for all $n_1\leq n$, then
\begin{align}
  F_{7}\left[A_2-nI\right]&=  F_{7}- x_2  B_2  \Big[\sum_{n_1=0}^{n-1}F_{7}\left[A_2-n_1 I, B_2+ I,  C_1+I\right]\Big] C^{-1}_{1}\notag\\& \quad -x_3   \Big[\sum_{n_1=0}^{n-1}F_{7}\left[ A_2-n_1 I,  B_3+ I, C_1+I\right]\Big]B_3 C^{-1}_{1}.
\label{c43eq90}\end{align}
\end{theorem}
\begin{theorem}
Let $A_2+nI$ be an invertible matrix for all $n\geq 0$ and let $A_i B_2= B_2 A_i$, $i=1, 2$; $B_1 B_2= B_2 B_1$;  $C_1 B_3= B_3 C_1$. Then the following recursion formula holds true for Lauricella  matrix function $F_7$:
\begin{align}
F_7\left[A_2+nI\right]&=\sum_{N_2\leq n}^{}{n\choose n_1, n_2} (B_2)_{n_1} x_{2}^{n_1}x_{3}^{n_2}\notag\\& \quad \times \Big[ F_{7}\left[A_2+N_2I, B_2+n_1I,  B_3+n_2I,  C_1+N_2 I\right]\Big](B_3)_{n_2}(C_1)^{-1}_{N_2}.
\label{c43eq91}\end{align}
Furthermore, if $A_2-n_1 I$ is invertible for all $n_1\leq n$, then
\begin{align}
F_7\left[A_2-nI\right]&=\sum_{N_2\leq n}^{}{n\choose n_1, n_2}  (B_2)_{n_1}  (-x_{2})^{n_1}(-x_{3})^{n_2}\notag\\& \quad \times \Big[ F_{7}\left[ B_2+n_1I,  B_3+n_2I,  C_1+N_2 I\right]\Big](B_3)_{n_2}(C_1)^{-1}_{N_2}.
\label{c43eq92}\end{align}
\end{theorem}
\begin{theorem}
Let $B_1+nI$ be an invertible matrix for all $n\geq 0$. Then the following recursion formula holds true for Lauricella  matrix function $F_7$:
\begin{align}
& F_{7}\left[B_1+nI\right] =  F_{7}+ x_1 A_1\Big[\sum_{n_1=1}^{n}F_{7}\left[A_1+ I, B_1+n_1 I,  C_1+I\right]\Big]C^{-1}_{1}.
\label{c43eq93}\end{align}
Furthermore, if $B_1-n_1 I$ is invertible for all $n_1\leq n$, then
\begin{align}
&F_{7}\left[B_1-nI\right]=  F_{7}- x_1 A_1\Big[\sum_{n_1=0}^{n-1}F_{7}\left[A_1 +I, B_1-n_1 I,  C_1+I\right]\Big]C^{-1}_{1}.
\label{c43eq94}\end{align}
\end{theorem}
\begin{theorem}
	Let $B_1+nI$ be an invertible matrix for all $n\geq0$. Then the following recursion formula holds true for Lauricella  matrix function $F_7$:
\begin{align}
& F_7\left[B_1+nI\right] =\sum_{n_1=0}^{n}{n\choose n_1} (A_1)_{n_1} x_{1}^{n_1} F_{7}\left[A_1+n_1I, B_1+n_1I, C_1+n_1I\right]\Big](C_1)^{-1}_{n_1}.
\label{c43eq95}\end{align}
Furthermore, if $B_1-n_1 I$ is invertible for all $n_1\leq n$, then
\begin{align}
& F_7\left[B_1-nI\right] =\sum_{n_1=0}^{n}{n\choose n_1} (A_1)_{n_1}(- x_{1})^{n_1} F_{7}\left[ A_1+n_1I, C_1+n_1I\right]\Big](C_1)^{-1}_{n_1}.
\label{c43eq96}\end{align}
\end{theorem}
\begin{theorem}
Let $B_i+nI$, $i= 2, 3$ be an invertible matrix for all $n\geq 0$ and let $A_1 A_2= A_2 A_1$. Then the following recursion formula holds true for Lauricella  matrix function $F_7$:
\begin{align}
& F_{7}\left[B_i+nI\right] =  F_{7}+ x_i A_2\Big[\sum_{n_1=1}^{n}F_{7}\left[A_2+ I, B_i+n_1 I,  C_1+I\right]\Big]C^{-1}_{1}.
\label{c43eq97}\end{align}
Furthermore, if $B_i-n_1 I$ is invertible for all $n_1\leq n$, then
\begin{align}
&F_{7}\left[B_i-nI\right]=  F_{7}- x_i A_2\Big[\sum_{n_1=0}^{n-1}F_{7}\left[A_2 +I, B_i-n_1 I,  C_1+I\right]\Big]C^{-1}_{1}.
\label{c43eq98}\end{align}
\end{theorem}
\begin{theorem}
Let $B_i+nI$,  $i= 2, 3$ be an invertible matrix for all $n\geq 0$ and let $A_1 A_2= A_2 A_1$. Then the following recursion formula holds true for Lauricella  matrix function $F_7$:
\begin{align}
& F_7\left[B_i+nI\right] =\sum_{n_1=0}^{n}{n\choose n_1} (A_2)_{n_1} x_{i}^{n_1} F_{7}\left[A_2+n_1I, B_i+n_1I, C_1+n_1I\right]\Big](C_1)^{-1}_{n_1}.
\label{c43eq99}\end{align}
Furthermore, if $B_i-n_1 I$ is invertible for all $n_1\leq n$, then
\begin{align}
& F_7\left[B_i-nI\right] =\sum_{n_1=0}^{n}{n\choose n_1} (A_2)_{n_1}(- x_{i})^{n_1} F_{7}\left[ A_2+n_1I, C_1+n_1I\right]\Big](C_1)^{-1}_{n_1}.
\label{c43eq100}\end{align}
\end{theorem}
\begin{theorem}
Let $C_1-nI$ be an invertible matrix for all $n\geq 0$ and let $A_1 A_2= A_2 A_1$, $A_i B_1= B_1 A_i$, $i=1, 2$; $B_j C_1= C_1 B_j$, $j=2, 3$. Then the following recursion formula holds true for Lauricella  matrix function $F_{7}$:
\begin{align}
\notag& F_{7}\left[C_1-nI\right]\\&= F_{7} +x_1 A_1 B_1\Big[\sum_{n_1=1}^{n}{F_7\left[A_1+I, B_1+I, C_1+(2-n_1)I\right]}{(C_1-n_1I)^{-1}(C_1-(n_1-1)I)^{-1}}\Big]\notag\\
& \quad +x_2 A_2 \Big[\sum_{n_1=1}^{n}{F_7 \left[A_2+I, B_2+I, C_1+(2-n_1)I\right]}B_2{(C_1-n_1I)^{-1}(C_1-(n_1-1)I)^{-1}}\Big]\notag\\
&\quad +x_3 A_2 \Big[\sum_{n_1=1}^{n}{F_7 \left[A_2+I, B_3+I, C_1+(2-n_1)I\right]}B_3{(C_1-n_1I)^{-1}(C_1-(n_1-1)I)^{-1}}\Big].\label{c43eq101}
\end{align}
\end{theorem}

\subsection*{ Recursion formulas of $F_8$:}

\begin{theorem}
Let $A_1+nI$ be an invertible matrix for all $n\geq 0$ and let $A_1 B_i =B_i  A_1$, $B_i B_j= B_j B_i$, $i, j=1,2,3$ and $C_1 C_2=C_2 C_1$. Then the following recursion formula holds true for Lauricella  matrix function $F_8$:
\begin{align}
&\notag F_8\left[{A_1}+nI\right]\\
&=F_8+ x_1 B_1 \left[\sum_{n_1=1}^{n}F_8\left[{A_1}+n_1 I, B_1+I, C_1+I\right]\right] C^{-1}_1\notag\\& \quad +  x_2 B_2 \left[\sum_{n_1=1}^{n}F_8\left[{A_1}+n_1 I, B_2+I, C_2+I\right]\right]C^{-1}_{2}\notag\\
& \quad +x_3 B_3\left[\sum_{n_1=1}^{n}F_8\left[{A_1}+n_1 I, B_3+I, C_2+I\right]\right]C^{-1}_2.\label{c43eq103}
\end{align}
Furthermore, if $A_1-n_1 I$ is invertible for all $n_1\leq n$, then
\begin{align}&
\notag F_8\left[{A_1}-n I\right]\\
&=F_8- x_1  B_{1}\left[\sum_{n_1=0}^{n-1}F_8\left[{A_1}-n_1 I, B_1+I,  C_1+I\right]\right]C^{-1}_1\notag\\& \quad - x_2 B_2\left[\sum_{n_1=0}^{n-1}F_8\left[{A_1}-n_1I, B_2+I, C_2+I\right]\right] C^{-1}_2\notag\\
& \quad -x_3  B_3\left[\sum_{n_1=0}^{n-1}F_8\left[{A_1}-n_1 I, B_3+I, C_2+I\right]\right] C^{-1}_2.\label{c43eq102}\end{align}
\end{theorem}
\begin{theorem}
Let $A_1+nI$ be an invertible matrix for all $n\geq 0$ and let $A_1 B_i =B_i  A_1$, $B_i B_j= B_j B_i$, $i, j=1,2,3$ and $C_1 C_2=C_2 C_1$. Then the following recursion formula holds true for  Lauricella  matrix function $F_8$:
\begin{align}
&F_8\left[A_1+nI\right] \nonumber\\
&=\sum_{N_3\leq n}^{}{n\choose n_1, n_2, n_3} \prod_{i=1}^{3}(B_i)_{n_i}\, x_{i}^{n_i}\notag\\&\Big[ F_{8}\left[A_1+N_3I, B_1+n_1I , B_2+n_2I, B_3+n_3I,  C_1+n_1I, C_2+(n_2+n_3) I\right]\Big](C_1)^{-1}_{n_1}(C_2)^{-1}_{n_2+n_3}.
\label{c43eq104}\end{align}
Furthermore, if $A_1-n_1 I$ is invertible for all $n_1\leq n$, then
\begin{align}
&F_8\left[A_1-nI\right]\nonumber\\
&=\sum_{N_3\leq n}^{}{n\choose n_1, n_2, n_3} \prod_{i=1}^{3}(B_i)_{n_i}\, (-x_{i})^{n_i}\notag\\&\Big[ F_{8}\left[ B_1+n_1I , B_2+n_2I, B_3+n_3I,   C_1+n_1I, C_2+(n_2+n_3) I\right]\Big](C_1)^{-1}_{n_1}(C_2)^{-1}_{n_2+n_3}.
\label{c43eq105}\end{align}
\end{theorem}
\begin{theorem}
Let $B_1+nI$ be an invertible matrix for all $n\geq 0$ and let $C_1C_2=C_2 C_1$. Then the following recursion formula holds true for Lauricella  matrix function $F_8$:
\begin{align}
& F_{8}\left[B_1+nI\right] =  F_{8} + x_1 A_1\Big[\sum_{n_1=1}^{n}F_{8}\left[A_1+ I, B_1+n_1 I,  C_1+I\right]\Big]C^{-1}_{1}.
\label{c43eq106}
\end{align}
Furthermore, if $B_1-n_1 I$ is invertible for all $n_1\leq n$, then
\begin{align}
&F_{8}\left[B_1-nI\right] =  F_{8}- x_1 A_1\Big[\sum_{n_1=0}^{n-1}F_{8}\left[A_1 +I, B_1-n_1 I,  C_1+I\right]\Big]C^{-1}_{1}.
\label{c43eq107}
\end{align}
\end{theorem}
\begin{theorem}
Let $B_1+nI$ be an invertible matrix for all $n\geq0$ and  let $C_1C_2=C_2 C_1$. Then the following recursion formula holds true for Lauricella  matrix function $F_8$:
\begin{align}
& F_8\left[B_1+nI\right] =\sum_{n_1=0}^{n}{n\choose n_1} (A_1)_{n_1} x_{1}^{n_1} F_{8}\left[A_1+n_1I, B_1+n_1I, C_1+n_1I\right]\Big](C_1)^{-1}_{n_1}.
\label{c43eq108}
\end{align}
Furthermore, if $B_1-n_1 I$ is invertible for all $n_1\leq n$, then
\begin{align}
& F_8\left[B_1-nI\right] =\sum_{n_1=0}^{n}{n\choose n_1} (A_1)_{n_1}(- x_{1})^{n_1} F_{8}\left[ A_1+n_1I, C_1+n_1I\right]\Big](C_1)^{-1}_{n_1}.
\label{c43eq109}
\end{align}
\end{theorem}
\begin{theorem}
Let $B_i+nI$, $i=2, 3$  be  invertible matrices for all $n\geq0$. Then the following recursion formula holds true for Lauricella  matrix function $F_8$:
\begin{align}
& F_{8}\left[B_i+nI\right] =  F_{8}+ x_i  A_1\Big[\sum_{n_1=1}^{n}F_{8}\left[A_1+ I, B_i+n_1 I,  C_2+I\right]\Big]C^{-1}_{2}.
\label{c43eq110}
\end{align}
Furthermore, if $B_i-n_1 I$ is invertible for all $n_1\leq n$, then
\begin{align}
&F_{8}\left[B_i-nI\right]=  F_{8}- x_i  A_1\Big[\sum_{n_1=0}^{n-1}F_{8}\left[A_1 +I, B_i-n_1 I,  C_2+I\right]\Big]C^{-1}_{2}.
\label{c43eq111}
\end{align}
\end{theorem}
\begin{theorem}
Let $B_i+nI$, $i=2, 3$ be invertible matrices for all $n\geq0$. Then the following recursion formula holds true for  Lauricella  matrix function $F_8$:
\begin{align}
& F_8\left[B_i+nI\right] =\sum_{n_1=0}^{n}{n\choose n_1} (A_1)_{n_1} x_{i}^{n_1} F_{8}\left[A_1+n_1I, B_i+n_1I, C_2+n_1I\right]\Big](C_2)^{-1}_{n_1}.
\label{c43eq112}
\end{align}
Furthermore, if $B_i-n_1 I$ is invertible for all $n_1\leq n$, then
\begin{align}
& F_8\left[B_i-nI\right] =\sum_{n_1=0}^{n}{n\choose n_1} (A_1)_{n_1}(- x_{i})^{n_1} F_{8}\left[ A_1+n_1I, C_2+n_1I\right]\Big](C_2)^{-1}_{n_1}.
\label{c43eq113}
\end{align}
\end{theorem}
\begin{theorem}
 Let $C_1-nI$   be an invertible matrix for all $n\geq0$ and let $A_1 B_1 =B_1  A_1$. Then the following recursion formula holds true for Lauricella  matrix function $F_{8}$:
\begin{align}
& F_{8}\left[{C_1}-nI\right] \notag\\& = F_{8} + x_1  A_1 B_1 \notag\\& \quad \times \Big[ \sum_{n_1=1}^{n}{ F_{8}\left[A_1+I,  B_1+I, {C_1}+(2-n_1)I\right] }{(C_1-n_1I)^{-1}(C_1-(n_1-1)I)^{-1}}\Big].
\label{c43eq114}
\end{align}
\end{theorem}
\begin{theorem}
Let $C_2-nI$ be an invertible matrix for all $n\geq 0$ and let $A_1 B_2= B_2 A_1$, $B_1 B_2= B_2 B_1$, $B_3 C_i= C_i B_3$, $i=1, 2$. Then the following recursion formula holds true for Lauricella  matrix function $F_{8}$:
\begin{align}
\notag& F_{8}\left[C_2-nI\right]\\&= F_{8} +x_2 A_1 B_2\Big[\sum_{n_1=1}^{n}{F_8\left[A_1+I, B_2+I, C_2+(2-n_1)I\right]}{(C_2-n_1I)^{-1}(C_2-(n_1-1)I)^{-1}}\Big]\notag\\
& \quad + x_3 A_1 \Big[\sum_{n_1=1}^{n}{F_8\left[A_1+I, B_3+I, C_2+(2-n_1)I\right]}B_3{(C_2-n_1I)^{-1}(C_2-(n_1-1)I)^{-1}}\Big]
.\label{c43eq115}
\end{align}
\end{theorem}

\subsection*{Recursion formulas of $F_{10}$:}

\begin{theorem}
Let $A_1+nI$ be an invertible matrix for all $n\geq 0$ and let $A_i B_1 = B_1 A_i$, $i=1, 2$; $C_1 C_2=C_2 C_1$. Then the following recursion formula holds true for Lauricella  matrix function $F_{10}$:
\begin{align}
 F_{10}\left[A_1+nI\right]&=  F_{10}+ x_1  B_1  \Big[\sum_{n_1=1}^{n}F_{10}\left[A_1+ n_1I, B_1+ I,  C_1+I\right]\Big] C^{-1}_{1}\notag\\& \quad + x_3 B_1  \Big[\sum_{n_1=1}^{n}F_{10}\left[  A_1+n_1 I, B_1+ I, C_2+I\right]\Big] C^{-1}_{2}.
\label{c43eq116}\end{align}
Furthermore, if $A_1-n_1 I$ is invertible for all $n_1\leq n$, then
\begin{align}
  F_{10}\left[A_1-nI\right]&=  F_{10}- x_1  B_1  \Big[\sum_{n_1=0}^{n-1}F_{10}\left[A_1-n_1 I, B_1+ I,  C_1+I\right]\Big] C^{-1}_{1}\notag\\& \quad - x_3 B_1  \Big[\sum_{n_1=0}^{n-1}F_{10}\left[ A_1-n_1 I,  B_1+ I, C_2+I\right]\Big] C^{-1}_{2}.
\label{c43eq117}
\end{align}
\end{theorem}
\begin{theorem}
Let $A_1+nI$ be an invertible matrix for all $n\geq 0$ and let $A_i B_1= B_1 A_i$, $i=1, 2$; $C_1 C_2=C_2 C_1$. Then the following recursion formula holds true for Lauricella  matrix function $F_{10}$:
\begin{align}
F_{10}\left[A_1+nI\right]&=\sum_{N_2\leq n}^{}{n\choose n_1, n_2} (B_1)_{N_2} x_{1}^{n_1}x_{3}^{n_2}\notag\\& \quad \times \Big[ F_{10}\left[A_1+N_2I, B_1+N_2I,   C_1+n_1 I, C_2+n_2I\right]\Big](C_1)^{-1}_{n_1}(C_2)^{-1}_{n_2}.
\label{c43eq118}
\end{align}
Furthermore, if $A_1-n_1 I$ is invertible for all $n_1\leq n$, then
\begin{align}
F_{10}\left[A_1-nI\right]&=\sum_{N_2\leq n}^{}{n\choose n_1, n_2}  (B_1)_{N_2}  (-x_{1})^{n_1}(-x_{3})^{n_2}\notag\\& \quad \times \Big[ F_{10}\left[ B_1+N_2I,  C_1+n_1 I, C_2+n_2I\right]\Big](C_1)^{-1}_{n_1}(C_2)^{-1}_{n_2}.
\label{c43eq119}\end{align}
\end{theorem}
Similarly,  recursion formulas  $F_{10}\left[{B_1}+n\right], F_{10}\left[{B_1}-n\right]$  for $F_{10}$ can be obtain by interchange of $A_1\leftrightarrow B_1$ in above theorems.
\begin{theorem}
Let $A_2+nI$ be an invertible matrix for all $n\geq0$ and let $C_i B_2=B_2 C_i$, $i=1, 2$. Then the following recursion formula holds true for Lauricella  matrix function $F_{10}$:
\begin{align}
& F_{10}\left[A_2+nI\right] =  F_{10} + x_2  \Big[\sum_{n_1=1}^{n}F_{10}\left[A_2+n_1 I, B_2+ I,  C_2+I\right]\Big]B_2 C^{-1}_{2}.
\label{c43eq120}\end{align}
Furthermore, if $A_2-n_1 I$ is invertible for all $n_1\leq n$, then
\begin{align}
&F_{10}\left[A_2-nI\right]=  F_{10}- x_2 \Big[\sum_{n_1=0}^{n-1}F_{10}\left[A_2 +n_1I, B_2+ I,  C_2+I\right]\Big]B_2C^{-1}_{2}.
\label{c43eq121}\end{align}
\end{theorem}
\begin{theorem}
Let $A_2+nI$ be an invertible matrix for all $n\geq 0$ and let $C_i B_2=B_2 C_i$, $i=1, 2$. Then the following recursion formula holds true for Lauricella  matrix function $F_{10}$:
\begin{align}
& F_{10}\left[A_2+nI\right] =\sum_{n_1=0}^{n}{n\choose n_1}  x_{2}^{n_1} F_{10}\left[A_2+n_1I, B_2+n_1I, C_2+n_1I\right]\Big](B_2)_{n_1}(C_2)^{-1}_{n_1}.
\label{c43eq122}\end{align}
Furthermore, if $A_2-n_1 I$ is invertible for all $n_1\leq n$, then
\begin{align}
& F_{10}\left[A_2-nI\right] = \sum_{n_1=0}^{n}{n\choose n_1} (- x_{2})^{n_1} F_{10}\left[ B_2+n_1I, C_2+n_1I\right]\Big](B_2)_{n_1}(C_2)^{-1}_{n_1}.
\label{c43eq123}
\end{align}
\end{theorem}
Similarly,  recursion formulas  $F_{10}\left[{B_2}+n\right], F_{10}\left[{B_2}-n\right]$  for $F_{10}$ can be obtained by interchange of $A_2\leftrightarrow B_2$ in above theorems.
\begin{theorem}
Let $C_1-nI$   be an invertible matrix for all $n\geq 0$ and let $A_i B_1 =B_1  A_i$, $i=1, 2$. Then the following recursion formula holds true for Lauricella  matrix function $F_{10}$:
\begin{align}
& F_{10}\left[{C_1}-nI\right] \notag\\&= F_{10} + x_1  A_1 B_1 \notag\\& \quad \times \Big[ \sum_{n_1=1}^{n}{ F_{10}\left[A_1+I,  B_1+I, {C_1}+(2-n_1)I\right] }{(C_1-n_1I)^{-1}(C_1-(n_1-1)I)^{-1}}\Big].
\label{c43eq124}\end{align}
\end{theorem}
\begin{theorem}
Let $C_2-nI$ be an invertible matrix for all $n\geq0$ and let $A_1 A_2= A_2 A_1$, $A_i B_1 =B_1  A_i$, $i=1, 2$ ; $B_2 C_j= C_j B_2$, $j=1, 2$. Then the following recursion formula holds true for Lauricella  matrix function $F_{10}$:
\begin{align}
\notag& F_{10}\left[C_2-nI\right]\\&= F_{10} +x_2 A_2\Big[\sum_{n_1=1}^{n}{F_{10}\left[A_2+I, B_2+I, C_2+(2-n_1)I\right]}B_2{(C_2-n_1I)^{-1}(C_2-(n_1-1)I)^{-1}}\Big]\notag\\
& \quad +x_3A_1  B_1 \Big[\sum_{n_1=1}^{n}{F_{10}\left[A_1+I, B_1+I, C_2+(2-n_1)I\right]}{(C_2-n_1I)^{-1}(C_2-(n_1-1)I)^{-1}}\Big]
.\label{c43eq125}
\end{align}
\end{theorem}

\subsection*{Recursion formulas of $F_{11}$:}

\begin{theorem}
Let $A_1+nI$ be an invertible matrix for all $n\geq 0$ and let $A_i B_1=B_1 A_i$, $i=1, 2$; $C_1C_2= C_2 C_1$. Then the following recursion formula holds true for Lauricella  matrix function $F_{11}$:
\begin{align}
& F_{11}\left[A_1+nI\right] =  F_{11}+ x_1 B_1 \Big[\sum_{n_1=1}^{n}F_{11}\left[A_1+n_1 I, B_1+ I,  C_1+I\right]\Big] C^{-1}_{1}.
\label{c43eq126}
\end{align}
Furthermore, if $A_1-n_1 I$ is invertible for all $n_1\leq n$, then
\begin{align}
&F_{11}\left[A_1-nI\right]=  F_{11}- x_1 B_1\Big[\sum_{n_1=0}^{n-1}F_{11}\left[A_1-n_1I, B_1+ I,  C_1+I\right]\Big]C^{-1}_{1}.
\label{c43eq127}\end{align}
\end{theorem}
\begin{theorem}
Let $A_1+nI$ be an invertible matrix for all $n\geq0$ and let $A_i B_1=B_1 A_i$, $i=1, 2$; $C_1C_2= C_2 C_1$. Then the following recursion formula holds true for Lauricella  matrix function $F_{11}$:
\begin{align}
& F_{11}\left[A_1+nI\right] =\sum_{n_1=0}^{n}{n\choose n_1}  x_{1}^{n_1} (B_1)_{n_1} F_{11}\left[A_1+n_1I, B_1+n_1I, C_1+n_1I\right]\Big](C_1)^{-1}_{n_1}.
\label{c43eq128}\end{align}
Furthermore, if $A_1-n_1 I$ is invertible for all $n_1\leq n$, then
\begin{align}
& F_{11}\left[A_1-nI\right] =\sum_{n_1=0}^{n}{n\choose n_1} (- x_{1})^{n_1} (B_1)_{n_1}  F_{11}\left[ B_1+n_1I, C_1+n_1I\right]\Big](C_1)^{-1}_{n_1}.
\label{c43eq129}\end{align}
\end{theorem}
\begin{theorem}
Let $A_2+nI$ be an invertible matrix for all $n\geq0$ and let $A_i B_1= B_1 A_i$, $C_i B_2=B_2 C_i$, $i=1, 2$. Then the following recursion formula holds true for Lauricella  matrix function $F_{11}$:
\begin{align}
 F_{11}\left[A_2+nI\right]&=  F_{11}+ x_2    \Big[\sum_{n_1=1}^{n}F_{11}\left[A_2+ n_1I, B_2+ I,  C_2+I\right]\Big]B_2 C^{-1}_{2}\notag\\& \quad + x_3 B_1  \Big[\sum_{n_1=1}^{n}F_{11}\left[  A_2+n_1 I, B_1+ I, C_2+I\right]\Big] C^{-1}_{2}.
\label{c43eq130}\end{align}
Furthermore, if $A_2-n_1 I$ is invertible for all $n_1\leq n$, then
\begin{align}
  F_{11}\left[A_2-nI\right]& =  F_{11}- x_2    \Big[\sum_{n_1=0}^{n-1}F_{11}\left[A_2-n_1 I, B_2+ I,  C_2+I\right]\Big]B_2 C^{-1}_{2}\notag\\& \quad - x_3 B_1  \Big[\sum_{n_1=0}^{n-1}F_{11}\left[ A_2-n_1 I,  B_1+ I, C_2+I\right]\Big] C^{-1}_{2}.
\label{c43eq131}\end{align}
\end{theorem}
\begin{theorem}
Let $A_2+nI$ be an invertible matrix for all $n\geq0$ and let $A_i B_1 = B_1 A_i$, $C_i B_2=B_2 C_i$, $i=1, 2$. Then the following recursion formula holds true for Lauricella  matrix function $F_{11}$:
\begin{align}
F_{11}\left[A_2+nI\right]&=\sum_{N_2\leq n}^{}{n\choose n_1, n_2} (B_1)_{n_2} x_{2}^{n_1}x_{3}^{n_2}\notag\\& \quad \times \Big[ F_{11}\left[A_2+N_2I, B_1+n_2I, B_2+n_1I,   C_2+N_2I\right]\Big](B_2)_{n_1}(C_2)^{-1}_{N_2}.
\label{c43eq132}\end{align}
Furthermore, if $A_2-n_1 I$ is invertible for all $n_1\leq n$, then
\begin{align}
F_{11}\left[A_2-nI\right]&=\sum_{N_2\leq n}^{}{n\choose n_1, n_2}  (B_1)_{n_2}  (-x_{2})^{n_1}(-x_{3})^{n_2}\notag\\& \quad \times \Big[ F_{11}\left[ B_1+n_2I, B_2+n_1I,   C_2+N_2I\right]\Big](B_2)_{n_1}(C_2)^{-1}_{N_2}.
\label{c43eq133}\end{align}
\end{theorem}
\begin{theorem}
Let $B_1+nI$ be an invertible matrix for all $n\geq0$ and let $A_1 A_2= A_2 A_1$, $C_1 C_2=C_2 C_1$. Then the following recursion formula holds true for Lauricella  matrix function $F_{11}$:
\begin{align}
 F_{11}\left[B_1+nI\right]&=  F_{11}+ x_1   A_1 \Big[\sum_{n_1=1}^{n}F_{11}\left[A_1+ I, B_1+n_1 I,  C_1+I\right]\Big] C^{-1}_{1}\notag\\& \quad +x_3 A_2  \Big[\sum_{n_1=1}^{n}F_{11}\left[  A_2+ I, B_1+n_1 I, C_2+I\right]\Big] C^{-1}_{2}.
\label{c43eq134}\end{align}
Furthermore, if $B_1-n_1 I$ is invertible for all $n_1\leq n$, then
\begin{align}
  F_{11}\left[B_1-nI\right]&=  F_{11}- x_1 A_1   \Big[\sum_{n_1=0}^{n-1}F_{11}\left[A_1+ I, B_1+n_1 I,  C_1+I\right]\Big] C^{-1}_{1}\notag\\& \quad - x_3 A_2  \Big[\sum_{n_1=0}^{n-1}F_{11}\left[ A_2+ I,  B_1-n_1 I, C_2+I\right]\Big] C^{-1}_{2}.
\label{c43eq135}\end{align}
\end{theorem}
\begin{theorem}
Let $B_1+nI$ be an invertible matrix for all $n\geq 0$ and let $A_1 A_2 = A_2 A_1$, $C_1 C_2 = C_2 C_1$. Then the following recursion formula holds true for Lauricella  matrix function $F_{11}$:
\begin{align}
F_{11}\left[B_1+nI\right]&=\sum_{N_2\leq n}^{}{n\choose n_1, n_2} (A_1)_{n_1}(A_2)_{n_2} x_{1}^{n_1}x_{3}^{n_2}\notag\\&\Big[ F_{11}\left[A_1+n_1I, A_+n_2I,  B_1+N_2I, C_1+n_1I,   C_2+n_2I\right]\Big](C_1)^{-1}_{n_1}(C_2)^{-1}_{n_2}.
\label{c43eq136}\end{align}
Furthermore, if $B_1-n_1 I$ is invertible for all $n_1\leq n$, then
\begin{align}
F_{11}\left[B_1-nI\right]&=\sum_{N_2\leq n}^{}{n\choose n_1, n_2} (A_1)_{n_1}(A_2)_{n_2} (-x_{1})^{n_1}(-x_{3})^{n_2}\notag\\& \quad \times \Big[ F_{11}\left[ A_1+n_1I, A_+n_2I,  C_1+n_1I,   C_2+n_2I\right]\Big](C_1)^{-1}_{n_1}(C_2)^{-1}_{n_2}.
\label{c43eq137}\end{align}
\end{theorem}
\begin{theorem}
Let $B_2+nI$ be an invertible matrix for all $n\geq0$ and let $A_1 A_2=A_2 A_1$. Then the following recursion formula holds true for Lauricella  matrix function $F_{11}$:
\begin{align}
& F_{11}\left[B_2+nI\right] =  F_{11}+ x_2 A_2 \Big[\sum_{n_1=1}^{n}F_{11}\left[A_2+ I, B_2+n_1 I,  C_2+I\right]\Big] C^{-1}_{2}.
\label{c43eq138}\end{align}
Furthermore, if $B_2-n_1 I$ is invertible for all $n_1\leq n$, then
\begin{align}
&F_{11}\left[B_2-nI\right]=  F_{11}- x_2 A_2\Big[\sum_{n_1=0}^{n-1}F_{11}\left[A_2+I, B_2-n_1 I,  C_2+I\right]\Big]C^{-1}_{2}.
\label{c43eq139}\end{align}
\end{theorem}
\begin{theorem}
Let $B_2+nI$ be an invertible matrix for all $n\geq0$ and let $A_1 A_2=A_2 A_1$. Then the following recursion formula holds true for Lauricella  matrix function $F_{11}$:
\begin{align}
& F_{11}\left[B_2+nI\right] =\sum_{n_1=0}^{n}{n\choose n_1}  x_{2}^{n_1} (A_2)_{n_1} F_{11}\left[A_2+n_1I, B_2+n_1I, C_2+n_1I\right]\Big](C_2)^{-1}_{n_1}.
\label{c43eq140}\end{align}
Furthermore, if $B_2-n_1 I$ is invertible for all $n_1\leq n$, then
\begin{align}
& F_{11}\left[B_2-nI\right] =\sum_{n_1=0}^{n}{n\choose n_1} (- x_{2})^{n_1} (A_2)_{n_1}  F_{11}\left[ A_2+n_1I, C_2+n_1I\right]\Big](C_2)^{-1}_{n_1}.
\label{c43eq141}\end{align}
\end{theorem}
\begin{theorem} Let $C_1-nI$   be an invertible matrix for all $n\geq0$ and let $A_i B_1 =B_1  A_i$, $i=1, 2$. Then the following recursion formula holds true for Lauricella  matrix function $F_{11}$:
\begin{align}
& F_{11}\left[{C_1}-nI\right] \notag\\& = F_{11} + x_1  A_1 B_1\notag\\& \quad \times \Big[ \sum_{n_1=1}^{n}{ F_{11}\left[A_1+I,  B_1+I, {C_1}+(2-n_1)I\right] }{(C_1-n_1I)^{-1}(C_1-(n_1-1)I)^{-1}}\Big].
\label{c43eq142}\end{align}
\end{theorem}
\begin{theorem}
Let $C_2-nI$ be an invertible matrix for all $n\geq 0$ and let $A_1 A_2 = A_2 A_1$, $A_i B_1 =B_1  A_i$, $B_2 C_i= C_i B_2$, $i=1, 2$. Then the following recursion formula holds true for  Lauricella  matrix function $F_{11}$:
\begin{align}
\notag& F_{11}\left[C_2-nI\right]\\&= F_{11} +x_2 A_2\Big[\sum_{n_1=1}^{n}{F_{11}\left[A_2+I, B_2+I, C_2+(2-n_1)I\right]}B_2{(C_2-n_1I)^{-1}(C_2-(n_1-1)I)^{-1}}\Big]\notag\\
& \quad + x_3A_1  B_1 \Big[\sum_{n_1=1}^{n}{F_{11}\left[A_1+I, B_1+I, C_2+(2-n_1)I\right]}{(C_2-n_1I)^{-1}(C_2-(n_1-1)I)^{-1}}\Big].\label{c43eq143}
\end{align}
\end{theorem}

\subsection*{Recursion formulas of $F_{12}$:}

\begin{theorem}
Let $A_1+nI$ be an invertible matrix for all $n\geq0$ and let $A_i B_1= B_1 A_i$, $C_i B_2=B_2 C_i$, $i=1, 2$. Then the following recursion formula holds true for Lauricella  matrix function $F_{12}$:
\begin{align}
 F_{12}\left[A_1+nI\right]&=  F_{12}+ x_1  B_1  \Big[\sum_{n_1=1}^{n}F_{12}\left[A_1+ n_1I, B_1+ I,  C_1+I\right]\Big] C^{-1}_{1}\notag\\& \quad + x_3    \Big[\sum_{n_1=1}^{n}F_{12}\left[  A_1+n_1 I, B_2+ I, C_2+I\right]\Big]B_2 C^{-1}_{2}.
\label{c43eq144}
\end{align}
Furthermore, if $A_1-n_1 I$ is invertible for all $n_1\leq n$, then
\begin{align}
  F_{12}\left[A_1-nI\right]&=  F_{12}- x_1 B_1    \Big[\sum_{n_1=0}^{n-1}F_{12}\left[A_1-n_1 I, B_1+ I,  C_1+I\right]\Big] C^{-1}_{1}\notag\\& \quad - x_3  \Big[\sum_{n_1=0}^{n-1}F_{12}\left[ A_1-n_1 I,  B_2+ I, C_2+I\right]\Big]B_2 C^{-1}_{2}.
\label{c43eq145}\end{align}
\end{theorem}
\begin{theorem}
Let $A_1+nI$ be an invertible matrix for all $n\geq0$ and let $A_i B_1= B_1 A_i$, $C_i B_2=B_2 C_i$, $i=1, 2$. Then the following recursion formula holds true for  Lauricella  matrix function $F_{12}$:
\begin{align}
&F_{12}\left[A_1+nI\right]\notag\\&=\sum_{N_2\leq n}^{}{n\choose n_1, n_2} (B_1)_{n_1} x_{1}^{n_1}x_{3}^{n_2}\notag\\& \quad \times \Big[ F_{12}\left[A_1+N_2I, B_1+n_1I, B_2+n_2I,   C_1+n_1I, C_2+n_2I\right]\Big](B_2)_{n_2}(C_1)^{-1}_{n_1}(C_2)^{-1}_{n_2}.\label{c43eq146}
\end{align}
Furthermore, if $A_1-n_1 I$ is invertible for all $n_1\leq n$, then
\begin{align}
F_{12}\left[A_1-nI\right]&=\sum_{N_2\leq n}^{}{n\choose n_1, n_2}  (B_1)_{n_1}  (-x_{1})^{n_1}(-x_{3})^{n_2}\notag\\&\Big[ F_{12}\left[ B_1+n_1I, B_2+n_2I,   C_1+n_1I, C_2+n_2I\right]\Big](B_2)_{n_2}(C_1)^{-1}_{n_1}(C_2)^{-1}_{n_2}.
\label{c43eq147}\end{align}
\end{theorem}
\begin{theorem}
Let $A_2+nI$ be an invertible matrix for all $n\geq0$ and let $A_i B_1=B_1 A_i$, $i=1, 2$. Then the following recursion formula holds true for Lauricella  matrix function $F_{12}$:
\begin{align}
& F_{12}\left[A_2+nI\right] =  F_{12}+ x_2 B_1 \Big[\sum_{n_1=1}^{n}F_{12}\left[A_1+n_1 I, B_1+ I,  C_2+I\right]\Big] C^{-1}_{2}.
\label{c43eq148}\end{align}
Furthermore, if $A_2-n_1 I$ is invertible for all $n_1\leq n$, then
\begin{align}
&F_{12}\left[A_2-nI\right]=  F_{12}- x_2 B_1\Big[\sum_{n_1=0}^{n-1}F_{12}\left[A_1-n_1I, B_1+ I,  C_2+I\right]\Big]C^{-1}_{2}.
\label{c43eq149}\end{align}
\end{theorem}
\begin{theorem}
Let $A_2+nI$ be an invertible matrix for all $n\geq0$ and let $A_i B_1=B_1 A_i$, $i=1, 2$. Then the following recursion formula holds true for Lauricella  matrix function $F_{12}$:
\begin{align}
& F_{12}\left[A_2+nI\right] =\sum_{n_1=0}^{n}{n\choose n_1}  x_{2}^{n_1} (B_1)_{n_1} F_{12}\left[A_1+n_1I, B_1+n_1I, C_2+n_1I\right]\Big](C_2)^{-1}_{n_1}.
\label{c43eq150}\end{align}
Furthermore, if $A_2-n_1 I$ is invertible for all $n_1\leq n$, then
\begin{align}
& F_{12}\left[A_2-nI\right] =\sum_{n_1=0}^{n}{n\choose n_1} (- x_{2})^{n_1} (B_1)_{n_1}  F_{12}\left[ B_1+n_1I, C_2+n_1I\right]\Big](C_2)^{-1}_{n_1}.
\label{c43eq151}\end{align}
\end{theorem}
\begin{theorem}
Let $B_1+nI$ be an invertible matrix for all $n\geq0$ and let $A_1A_2= A_2 A_1$, $C_1 C_2=C_2 C_1$. Then the following recursion formula holds true for Lauricella  matrix function $F_{12}$:
\begin{align}
 F_{12}\left[B_1+nI\right]&=  F_{12}+ x_1  A_1  \Big[\sum_{n_1=1}^{n}F_{12}\left[A_1+ I, B_1+n_1 I,  C_1+I\right]\Big] C^{-1}_{1}\notag\\& \quad +x_2  A_2   \Big[\sum_{n_1=1}^{n}F_{12}\left[  A_2+I, B_1+n_1 I, C_2+I\right]\Big] C^{-1}_{2}.
\label{c43eq152}\end{align}
Furthermore, if $B_1-n_1 I$ is invertible for all $n_1\leq n$, then
\begin{align}
  F_{12}\left[B_1-nI\right]&=  F_{12}- x_1 A_1    \Big[\sum_{n_1=0}^{n-1}F_{12}\left[A_1+ I, B_1+ I,  C_1+I\right]\Big] C^{-1}_{1}\notag\\& \quad -x_2 A_2 \Big[\sum_{n_1=0}^{n-1}F_{12}\left[ A_2+ I,  B_1-n_1 I, C_2+I\right]\Big] C^{-1}_{2}.
\label{c43eq153}\end{align}
\end{theorem}
\begin{theorem}
Let $B_1+nI$ be an invertible matrix for all $n\geq0$ and let $A_1A_2= A_2 A_1$, $C_1 C_2=C_2 C_1$. Then the following recursion formula holds true for Lauricella  matrix function $F_{12}$:
\begin{align}
&F_{12}\left[B_1+nI\right]\nonumber\\
&=\sum_{N_2\leq n}^{}{n\choose n_1, n_2} (A_1)_{n_1}(A_2)_{n_2}\, x_{1}^{n_1}x_{2}^{n_2}\notag\\& \quad \times \Big[ F_{12}\left[A_1+n_1I, A_2+n_2I,  B_1+N_2I,   C_1+n_1I, C_2+n_2I\right]\Big](C_1)^{-1}_{n_1}(C_2)^{-1}_{n_2}.
\label{c43eq154}\end{align}
Furthermore, if $B_1-n_1 I$ is invertible for all $n_1\leq n$, then
\begin{align}
F_{12}\left[B_1-nI\right]&=\sum_{N_2\leq n}^{}{n\choose n_1, n_2} (A_1)_{n_1}(A_2)_{n_2}  (-x_{1})^{n_1}(-x_{2})^{n_2}\notag\\& \quad \times \Big[ F_{12}\left[ A_1+n_1I, A_2+n_2I,   C_1+n_1I, C_2+n_2I\right]\Big](C_1)^{-1}_{n_1}(C_2)^{-1}_{n_2}.
\label{c43eq155}\end{align}
\end{theorem}
\begin{theorem}
Let $B_2+nI$ be an invertible matrix for all $n\geq0$. Then the following recursion formula holds true for Lauricella  matrix function $F_{12}$:
\begin{align}
& F_{12}\left[B_2+nI\right] =  F_{12}+ x_3 A_1 \Big[\sum_{n_1=1}^{n}F_{12}\left[A_1+ I, B_2+n_1 I,  C_2+I\right]\Big] C^{-1}_{2}.
\label{c43eq156}\end{align}
Furthermore, if $B_2-n_1 I$ is invertible for all $n_1\leq n$, then
\begin{align}
&F_{12}\left[B_2-nI\right]=  F_{12}- x_3 A_1\Big[\sum_{n_1=0}^{n-1}F_{12}\left[A_1+I, B_2-n_1 I,  C_2+I\right]\Big]C^{-1}_{2}.
\label{c43eq157}\end{align}
\end{theorem}
\begin{theorem}
Let $B_2+nI$ be an invertible matrix for all $n\geq0$. Then the following recursion formula holds true for generalized Lauricella  matrix function $F_{12}$:
\begin{align}
& F_{12}\left[B_2+nI\right] =\sum_{n_1=0}^{n}{n\choose n_1}  x_{3}^{n_1} (A_1)_{n_1} F_{12}\left[A_1+n_1I, B_2+n_1I, C_2+n_1I\right]\Big](C_2)^{-1}_{n_1}.
\label{c43eq158}\end{align}
Furthermore, if $B_2-n_1 I$ is invertible for all $n_1\leq n$, then
\begin{align}
& F_{12}\left[B_2-nI\right] =\sum_{n_1=0}^{n}{n\choose n_1} (- x_{3})^{n_1} (A_1)_{n_1}  F_{12}\left[ A_1+n_1I, C_2+n_1I\right]\Big](C_2)^{-1}_{n_1}.
\label{c43eq159}\end{align}
\end{theorem}
\begin{theorem}
 Let $C_1-nI$   be an invertible matrix for all $n\geq0$ and let $A_i B_1 =B_1  A_i$, $i=1, 2$, $C_1 C_2=C_2 C_1$.Then the following recursion formula holds true for Lauricella  matrix function $F_{12}$:
\begin{align}
& F_{12}\left[{C_1}-nI\right] \notag\\& = F_{12} + x_1  A_1 B_1\notag\\& \quad \times \Big[ \sum_{n_1=1}^{n}{ F_{12}\left[A_1+I,  B_1+I, {C_1}+(2-n_1)I\right] }{(C_1-n_1I)^{-1}(C_1-(n_1-1)I)^{-1}}\Big].
\label{c43eq160}\end{align}
\end{theorem}
\begin{theorem}
Let $C_2-nI$ be an invertible matrix for all $n\geq 0$ and let $A_1 A_2 = A_2 A_1$; $A_i B_1 =B_1  A_i$, $B_2 C_i= C_i B_2$, $i=1, 2$. Then the following recursion formula holds true for Lauricella  matrix function $F_{12}$:
\begin{align}
\notag& F_{12}\left[C_2-nI\right]\\&= F_{12} +x_2 A_2B_1\Big[\sum_{n_1=1}^{n}{F_{12}\left[A_2+I, B_1+I, C_2+(2-n_1)I\right]}{(C_2-n_1I)^{-1}(C_2-(n_1-1)I)^{-1}}\Big]\notag\\
& \quad + x_3A_1  \Big[\sum_{n_1=1}^{n}{F_{12}\left[A_1+I, B_2+I, C_2+(2-n_1)I\right]}B_2{(C_2-n_1I)^{-1}(C_2-(n_1-1)I)^{-1}}\Big]
.\label{c43eq161}
\end{align}
\end{theorem}

\subsection*{Recursion formulas of $F_{13}$:}

\begin{theorem}
Let $A_1+nI$ be an invertible matrix for all $n\geq0$ and let $A_i B_1=B_1 A_i$, $i=1, 2$. Then the following recursion formula holds true for Lauricella  matrix function $F_{13}$:
\begin{align}
& F_{13}\left[A_1+nI\right] =  F_{13}+ x_1 B_1 \Big[\sum_{n_1=1}^{n}F_{13}\left[A_1+n_1 I, B_1+ I,  C_1+I\right]\Big] C^{-1}_{1}.
\label{c43eq162}\end{align}
Furthermore, if $A_1-n_1 I$ is invertible for all $n_1\leq n$, then
\begin{align}
&F_{13}\left[A_1-nI\right] =  F_{13}- x_1 B_1\Big[\sum_{n_1=0}^{n-1}F_{13}\left[A_1-n_1I, B_1+ I,  C_1+I\right]\Big]C^{-1}_{1}.
\label{c43eq163}\end{align}
\end{theorem}
\begin{theorem}
Let $A_1+nI$ be an invertible matrix for all $n\geq0$ and let $A_i B_1=B_1 A_i$, $i=1, 2$. Then the following recursion formula holds true for Lauricella  matrix function $F_{13}$:
\begin{align}
& F_{13}\left[A_1+nI\right] =\sum_{n_1=0}^{n}{n\choose n_1}  x_{1}^{n_1} (B_1)_{n_1} F_{13}\left[A_1+n_1I, B_1+n_1I, C_1+n_1I\right]\Big](C_1)^{-1}_{n_1}.
\label{c43eq164}\end{align}
Furthermore, if $A_1-n_1 I$ is invertible for all $n_1\leq n$, then
\begin{align}
& F_{13}\left[A_1-nI\right] =\sum_{n_1=0}^{n}{n\choose n_1} (- x_{1})^{n_1} (B_1)_{n_1}  F_{13}\left[ B_1+n_1I, C_1+n_1I\right]\Big](C_1)^{-1}_{n_1}.
\label{c43eq165}\end{align}
\end{theorem}
\begin{theorem}
Let $A_2+nI$ be an invertible matrix for all $n\geq0$ and let $A_i B_1= B_1 A_i$, $i=1, 2$, $C_1 B_2=B_2 C_1$. Then the following recursion formula holds true for Lauricella  matrix function $F_{13}$:
\begin{align}
 F_{13}\left[A_2+nI\right]&=  F_{13}+ x_2   \Big[\sum_{n_1=1}^{n}F_{13}\left[A_2+ n_1I, B_2+ I,  C_1+I\right]\Big] B_2C^{-1}_{1}\notag\\& \quad +x_3   B_1 \Big[\sum_{n_1=1}^{n}F_{13}\left[  A_2+n_1 I, B_1+ I, C_1+I\right]\Big] C^{-1}_{1}.
\label{c43eq166}\end{align}
Furthermore, if $A_2-n_1 I$ is invertible for all $n_1\leq n$, then
\begin{align}
  F_{13}\left[A_2-nI\right]&=  F_{13}- x_2     \Big[\sum_{n_1=0}^{n-1}F_{13}\left[A_2-n_1 I, B_2+ I,  C_1+I\right]\Big]B_2 C^{-1}_{1}\notag\\&\quad -x_3 B_1 \Big[\sum_{n_1=0}^{n-1}F_{13}\left[ A_2-n_1 I,  B_1+ I, C_1+I\right]\Big]C^{-1}_{1}.
\label{c43eq167}\end{align}
\end{theorem}
\begin{theorem}
Let $A_2+nI$ be an invertible matrix for all $n\geq0$ and let $A_i B_1= B_1 A_i$, $i=1, 2$; $C_1 B_2=B_2 C_1$. Then the following recursion formula holds true for Lauricella  matrix function $F_{13}$:
\begin{align}
F_{13}\left[A_2+nI\right]&=\sum_{N_2\leq n}^{}{n\choose n_1, n_2} (B_1)_{n_2} x_{2}^{n_1}x_{3}^{n_2}\notag\\& \quad \times \Big[ F_{13}\left[A_2+N_2I, B_1+n_2I, B_2+n_1I,   C_1+N_2I\right]\Big](B_2)_{n_1}(C_1)^{-1}_{N_2}.
\label{c43eq168}\end{align}
Furthermore, if $A_2-n_1 I$ is invertible for all $n_1\leq n$, then
\begin{align}
F_{13}\left[A_2-nI\right]&=\sum_{N_2\leq n}^{}{n\choose n_1, n_2}  (B_1)_{n_2}  (-x_{2})^{n_1}(-x_{3})^{n_2}\notag\\& \quad \times \Big[ F_{13}\left[ B_1+n_2I, B_2+n_1I,   C_1+N_2I\right]\Big](B_2)_{n_1}(C_1)^{-1}_{N_2}.
\label{c43eq169}\end{align}
\end{theorem}
\begin{theorem}
Let $B_1+nI$ be an invertible matrix for all $n\geq0$ and let $A_1 A_2= A_2 A_1$. Then the following recursion formula holds true for Lauricella  matrix function $F_{13}$:
\begin{align}
 F_{13}\left[B_1+nI\right]&=  F_{13}+ x_1 A_1  \Big[\sum_{n_1=1}^{n}F_{13}\left[A_1+ I, B_1+n_1 I,  C_1+I\right]\Big] C^{-1}_{1}\notag\\& \quad +x_3   A_2 \Big[\sum_{n_1=1}^{n}F_{13}\left[  A_2+1 , B_1+n_1 I, C_1+I\right]\Big] C^{-1}_{1}.
\label{c43eq170}\end{align}
Furthermore, if $B_1-n_1 I$ is invertible for all $n_1\leq n$, then
\begin{align}
  F_{13}\left[B_1-nI\right]&=  F_{13}- x_1 A_1    \Big[\sum_{n_1=0}^{n-1}F_{13}\left[A_1+ I, B_1-n_1 I,  C_1+I\right]\Big] C^{-1}_{1}\notag\\& \quad -x_3 A_2 \Big[\sum_{n_1=0}^{n-1}F_{13}\left[ A_2+ I,  B_1-n_1 I, C_1+I\right]\Big]C^{-1}_{1}.
\label{c43eq171}\end{align}
\end{theorem}
\begin{theorem}
Let $B_1+nI$ be an invertible matrix for all $n\geq0$ and let $A_1 A_2= A_2 A_1$. Then the following recursion formula holds true for Lauricella  matrix function $F_{13}$:
\begin{align}
F_{13}\left[B_1+nI\right]&=\sum_{N_2\leq n}^{}{n\choose n_1, n_2} (A_1)_{n_1}(A_2)_{n_2} x_{1}^{n_1}x_{3}^{n_2}\notag\\& \quad \times \Big[ F_{13}\left[A_1+n_1I, A_2+n_2I,  B_1+N_2I,   C_1+N_2I\right]\Big](C_1)^{-1}_{N_2}.
\label{c43eq172}\end{align}
Furthermore, if $B_1-n_1 I$ is invertible for all $n_1\leq n$, then
\begin{align}
F_{13}\left[B_1-nI\right]&=\sum_{N_2\leq n}^{}{n\choose n_1, n_2}  (A_1)_{n_1}(A_2)_{n_2}  (-x_{1})^{n_1}(-x_{3})^{n_2}\notag\\& \quad \times \Big[ F_{13}\left[ A_1+n_1I, A_2+n_2I,    C_1+N_2I\right]\Big](C_1)^{-1}_{N_2}.
\label{c43eq173}\end{align}
\end{theorem}
\begin{theorem}
Let $B_2+nI$ be an invertible matrix for all $n\geq0$ and let $A_1 A_2=A_2 A_1$. Then the following recursion formula holds true for Lauricella  matrix function $F_{13}$:
\begin{align}
& F_{13}\left[B_2+nI\right] =  F_{13}+ x_2 A_2 \Big[\sum_{n_1=1}^{n}F_{13}\left[A_2+ I, B_2+n_1 I,  C_1+I\right]\Big] C^{-1}_{1}.
\label{c43eq174}\end{align}
Furthermore, if $B_2-n_1 I$ is invertible for all $n_1\leq n$, then
\begin{align}
&F_{13}\left[B_2-nI\right]=  F_{13}- x_2 A_2\Big[\sum_{n_1=0}^{n-1}F_{13}\left[A_2+I, B_2-n_1 I,  C_1+I\right]\Big]C^{-1}_{1}.
\label{c43eq175}\end{align}
\end{theorem}
\begin{theorem}
Let $B_2+nI$ be an invertible matrix for all $n\geq0$ and let $A_1 A_2=A_2 A_1$. Then the following recursion formula holds true for Lauricella  matrix function $F_{13}$:
\begin{align}
& F_{13}\left[B_2+nI\right] =\sum_{n_1=0}^{n}{n\choose n_1}  x_{2}^{n_1} (A_2)_{n_1} F_{13}\left[A_2+n_1I, B_2+n_1I, C_1+n_1I\right]\Big](C_1)^{-1}_{n_1}.
\label{c43eq176}\end{align}
Furthermore, if $B_2-n_1 I$ is invertible for all $n_1\leq n$, then
\begin{align}
& F_{13}\left[B_2-nI\right] =\sum_{n_1=0}^{n}{n\choose n_1} (- x_{2})^{n_1} (A_2)_{n_1}  F_{13}\left[ A_2+n_1I, C_1+n_1I\right]\Big](C_1)^{-1}_{n_1}.
\label{c43eq177}\end{align}
\end{theorem}
\begin{theorem}
Let $C_1-nI$ be an invertible matrix for all $n\geq0$ and let $A_1 A_2= A_2 A_1$, $A_i B_1 =B_1  A_i$, $i=1, 2$; $B_2 C_1= C_1 B_2$. Then the following recursion formula holds true for Lauricella  matrix function $F_{13}$:
\begin{align}
\notag& F_{13}\left[C_1-nI\right]\\&= F_{13} +x_1 A_1B_1\Big[\sum_{n_1=1}^{n}{F_{13}\left[A_1+I, B_1+I, C_1+(2-n_1)I\right]}{(C_1-n_1I)^{-1}(C_1-(n_1-1)I)^{-1}}\Big]\notag\\
& \quad +x_2 A_2 \Big[\sum_{n_1=1}^{n}{F_{13}\left[A_2+I, B_2+I, C_1+(2-n_1)I\right]}B_2{(C_1-n_1I)^{-1}(C_1-(n_1-1)I)^{-1}}\Big]\notag\\& \quad +x_3 A_2 B_1\Big[\sum_{n_1=1}^{n}{F_{13}\left[A_2+I, B_1+I, C_1+(2-n_1)I\right]}{(C_1-n_1I)^{-1}(C_1-(n_1-1)I)^{-1}}\Big]
.\label{c43eq178}
\end{align}
\end{theorem}

\subsection*{Recursion formulas of $F_{14}$:}

\begin{theorem}
Let $A_1+nI$ be an invertible matrix for all $n\geq0$ and let $A_1 B_1=B_1  A_1$,  $B_2 C_i= C_i B_2$, $i=1,2$; $C_1 C_2=C_2 C_1$. Then the following recursion formula holds true for Lauricella  matrix function $F_{14}$:
\begin{align}
&\notag F_{14}\left[{A_1}+nI\right]\\
&=F_{14}+ x_1 B_1 \left[\sum_{n_1=1}^{n}F_{14}\left[{A_1}+n_1 I, B_1+I, C_1+I\right]\right] C^{-1}_1\notag\\& \quad +  x_2  \left[\sum_{n_1=1}^{n}F_{14}\left[{A_1}+n_1 I, B_2+I, C_2+I\right]\right]B_2C^{-1}_{2}\notag\\
& \quad +x_3 B_1\left[\sum_{n_1=1}^{n}F_{14}\left[{A_1}+n_1 I, B_1+I, C_2+I\right]\right]C^{-1}_2.\label{c43eq179}
\end{align}
Furthermore, if $A_1-n_1 I$ is invertible for all $n_1\leq n$, then
\begin{align}&
\notag F_{14}\left[{A_1}-n I\right]\\
&=F_{14}- x_1  B_{1}\left[\sum_{n_1=0}^{n-1}F_{14}\left[{A_1}-n_1 I, B_1+I,  C_1+I\right]\right]C^{-1}_1\notag\\& \quad - x_2 \left[\sum_{n_1=0}^{n-1}F_{14}\left[{A_1}-n_1I, B_2+I, C_2+I\right]\right]B_2 C^{-1}_2\notag\\
& \quad -x_3  B_1\left[\sum_{n_1=0}^{n-1}F_{14}\left[{A_1}-n_1 I, B_1+I, C_2+I\right]\right] C^{-1}_2.\label{c43eq180}\end{align}
\end{theorem}
\begin{theorem}
Let $A_1+nI$ be an invertible matrix for all $n\geq0$ and let $A_1 B_1=B_1  A_1$,  $B_2 C_i= C_i B_2$, $i=1,2$; $C_1 C_2=C_2 C_1$. Then the following recursion formula holds true for Lauricella  matrix function $F_{14}$:
\begin{align}
F_{14}\left[A_1+nI\right]&=\sum_{N_3\leq n}^{}{n\choose n_1, n_2, n_3} (B_1)_{n_1+n_3}\prod_{i=1}^{3}\, x_{i}^{n_i}\notag\\& \quad \times \Big[ F_{{14}}\left[A_1+N_3I, B_1+(n_1+n_3)I , B_2+n_2I,  C_1+n_1I, C_2+(n_2+n_3) I\right]\Big]\notag\\& \quad \times (B_2)_{n_2}(C_1)^{-1}_{n_1}(C_2)^{-1}_{n_2+n_3}.
\label{c43eq181}\end{align}
Furthermore, if $A_1-n_1 I$ is invertible for all $n_1\leq n$, then
\begin{align}
F_{14}\left[A_1-nI\right]&=\sum_{N_3\leq n}^{}{n\choose n_1, n_2, n_3}(B_1)_{n_1+n_3} \prod_{i=1}^{3} (-x_{i})^{n_i}\notag\\& \quad \times \Big[ F_{{14}}\left[ B_1+(n_1+n_3)I , B_2+n_2I,  C_1+n_1I, C_2+(n_2+n_3) I\right]\Big]\notag\\& \quad \times (B_2)_{n_2}(C_1)^{-1}_{n_1}(C_2)^{-1}_{n_2+n_3}.
\label{c43eq182}\end{align}
\end{theorem}
\begin{theorem}
Let $B_1+nI$ be an invertible matrix for all $n\geq0$ and let $C_1 C_2= C_2 C_1$. Then the following recursion formula holds true for Lauricella  matrix function $F_{14}$:
\begin{align}
 F_{14}\left[B_1+nI\right]&=  F_{14}+ x_1 A_1  \Big[\sum_{n_1=1}^{n}F_{14}\left[A_1+ I, B_1+n_1 I,  C_1+I\right]\Big] C^{-1}_{1}\notag\\& \quad +x_3   A_1 \Big[\sum_{n_1=1}^{n}F_{14}\left[  A_1+I , B_1+n_1 I, C_2+I\right]\Big] C^{-1}_{2}.
\label{c43eq183}\end{align}
Furthermore, if $B_1-n_1 I$ is invertible for all $n_1\leq n$, then
\begin{align}
  F_{14}\left[B_1-nI\right]&=  F_{14}- x_1 A_1    \Big[\sum_{n_1=0}^{n-1}F_{14}\left[A_1+ I, B_1-n_1 I,  C_1+I\right]\Big] C^{-1}_{1}\notag\\& \quad -x_3 A_1 \Big[\sum_{n_1=0}^{n-1}F_{14}\left[ A_1+ I,  B_1-n_1 I, C_2+I\right]\Big]C^{-1}_{2}.
\label{c43eq184}\end{align}
\end{theorem}
\begin{theorem}
Let $B_1+nI$ be an invertible matrix for all $n\geq0$ and let $C_1 C_2= C_2 C_1$. Then the following recursion formula holds true for Lauricella  matrix function $F_{14}$:
\begin{align}
F_{14}\left[B_1+nI\right]&=\sum_{N_2\leq n}^{}{n\choose n_1, n_2} (A_1)_{N_2} x_{1}^{n_1}x_{3}^{n_2}\notag\\& \quad \times \Big[ F_{14}\left[A_1+N_2I,  B_1+N_2I,   C_1+n_1I, C_2+n_2I\right]\Big](C_1)^{-1}_{n_1}(C_2)^{-1}_{n_2}.
\label{c43eq185}\end{align}
Furthermore, if $B_1-n_1 I$ is invertible for all $n_1\leq n$, then
\begin{align}
F_{14}\left[B_1-nI\right]&=\sum_{N_2\leq n}^{}{n\choose n_1, n_2}  (A_1)_{N_2}(-x_{1})^{n_1}(-x_{3})^{n_2}\notag\\& \quad \times \Big[ F_{14}\left[A_1+N_2I,  C_1+n_1I, C_2+n_2I\right]\Big](C_1)^{-1}_{n_1}(C_2)^{-1}_{n_2}.
\label{c43eq186}\end{align}
\end{theorem}
\begin{theorem}
Let $B_2+nI$ be an invertible matrix for all $n\geq0$ . Then the following recursion formula holds true for Lauricella  matrix function $F_{14}$:
\begin{align}
& F_{14}\left[B_2+nI\right] =  F_{14}+ x_2 A_1 \Big[\sum_{n_1=1}^{n}F_{14}\left[A_1+ I, B_2+n_1 I,  C_2+I\right]\Big] C^{-1}_{2}.
\label{c43eq187}\end{align}
Furthermore, if $B_2-n_1 I$ is invertible for  all $n_1\leq n$, then
\begin{align}
&F_{14}\left[B_2-nI\right] =  F_{14}- x_2 A_1\Big[\sum_{n_1=0}^{n-1}F_{14}\left[A_1+I, B_2-n_1 I,  C_2+I\right]\Big]C^{-1}_{2}.
\label{c43eq188}\end{align}
\end{theorem}
\begin{theorem}
Let $B_2+nI$ be an invertible matrix for all $n\geq0$. Then the following recursion formula holds true for Lauricella  matrix function $F_{14}$:
\begin{align}
& F_{14}\left[B_2+nI\right] =\sum_{n_1=0}^{n}{n\choose n_1}  x_{2}^{n_1} (A_1)_{n_1} F_{14}\left[A_1+n_1I, B_2+n_1I, C_2+n_1I\right]\Big](C_2)^{-1}_{n_1}.
\label{c43eq189}\end{align}
Furthermore, if $B_2-n_1 I$ is invertible for all $n_1\leq n$, then
\begin{align}
& F_{14}\left[B_2-nI\right] =\sum_{n_1=0}^{n}{n\choose n_1} (- x_{2})^{n_1} (A_1)_{n_1}  F_{14}\left[ A_1+n_1I, C_2+n_1I\right]\Big](C_2)^{-1}_{n_1}.
\label{c43eq190}\end{align}
\end{theorem}
\begin{theorem}
 Let $C_1-nI$   be an invertible matrix for all $n\geq 0$ and let $A_1 B_1= B_1 A_1$. Then the following recursion formula holds true for Lauricella  matrix function $F_{14}$:
\begin{align}
& F_{14}\left[{C_1}-nI\right] \nonumber\\
& = F_{14} + x_1  A_1 B_1\nonumber\\
& \quad \times \Big[ \sum_{n_1=1}^{n}{ F_{14}\left[A_1+I,  B_1+I, {C_1}+(2-n_1)I\right] }{(C_1-n_1I)^{-1}(C_1-(n_1-1)I)^{-1}}\Big].
\label{c43eq191}\end{align}
\end{theorem}
\begin{theorem}Let $C_2-nI$ be an invertible matrix for all $n\geq0$ and let $A_1 B_1= B_1 A_1$, $B_2 C_i=C_i B_2$, $i=1, 2$. Then the following recursion formula holds true for Lauricella  matrix function $F_{14}$:
\begin{align}
\notag& F_{14}\left[C_2-nI\right]\\&= F_{14} +x_2 A_1\Big[\sum_{n_1=1}^{n}{F_{14}\left[A_1+I, B_2+I, C_2+(2-n_1)I\right]}B_2{(C_2-n_1I)^{-1}(C_2-(n_1-1)I)^{-1}}\Big]\notag\\
& \quad +x_3A_1 B_1 \Big[\sum_{n_1=1}^{n}{F_{14}\left[A_1+I, B_1+I, C_2+(2-n_1)I\right]}{(C_2-n_1I)^{-1}(C_2-(n_1-1)I)^{-1}}\Big].\label{c43eq192}
\end{align}
\end{theorem}

\end{document}